\documentclass{cmslatex}
\usepackage{latexsym, amssymb, enumerate, amsmath}
\usepackage{graphicx}
\usepackage{subfigure}
\usepackage{float}
\usepackage{epsf}
\usepackage{array}
\usepackage[procent]{overpic}
\usepackage[mathlines]{lineno}
\usepackage[top=1in, bottom=1in, left=1in, right=1in]{geometry}
\usepackage{setspace}
\usepackage{cite}
\usepackage{color}

\DeclareGraphicsExtensions{.pdf, .pdftex, .eps, .jpg, .png, {}}
\renewcommand{\theequation}{\arabic{section}.\arabic{equation}}

\def\dx{\Delta x}
\def\dt{\Delta t}
\def\xs{x^{\star}}
\def\dxs{\Delta\xs}
\def\jpmh{{j\pm\frac{1}{2}}}
\def\jph{{j+\frac{1}{2}}}
\def\jmh{{j-\frac{1}{2}}}
\def\CFL{\mathbf{CFL}}
\newcommand{\mU}{\mathbf{U}}

\newcommand{\mF}{\mathbf{F}}
\newcommand{\mS}{\mathbf{S}}
\newcommand{\mH}{\mathbf{H}}
\newcommand{\mo}{\mathbf{0}}
\newcommand{\sgn}{\mathrm{sgn}}

\newcommand\eref[1]{(\ref{#1})}

\newcommand{\dtd}{\dt_j^{\text{drain}}}

\newtheorem{thm}{Theorem}[section]

\newtheorem{lem}[thm]{Lemma}

\newtheorem{defn}[thm]{Definition}

\newtheorem{rem}[thm]{Remark}
\title{A well-balanced reconstruction of wet/dry fronts for the shallow water equations
}

\author{Andreas Bollermann \footnotemark[2]
\and    Guoxian Chen \footnotemark[3] \footnotemark[2]
\and    Alexander Kurganov \footnotemark[4]
\and    Sebastian Noelle \footnotemark[2]\footnotemark[1] }

\footnotetext[1]{Corresponding author. noelle@igpm.rwth-aachen.de}
\footnotetext[2]{IGPM, RWTH Aachen University of Technology, Templergraben 55, 52062, Aachen, Germany.}
\footnotetext[3]{School of Mathematics and Statistics, Wuhan University, Wuhan, 430072, P.R. China.}
\footnotetext[4]{Mathematics Department, Tulane University, New Orleans, LA 70118, USA.}

\begin{document}
\large
\begin{spacing}{1.2}
\maketitle
\begin{abstract}

In this paper, we construct a well-balanced, positivity preserving finite volume
scheme for the shallow water equations based on a continuous, piecewise linear
discretization of the bottom topography.  The main new technique is a special
reconstruction of the flow variables in wet-dry cells, which is presented in
this paper for the one dimensional case. We realize the new reconstruction in
the framework of the second-order semi-discrete central-upwind scheme from (A.
Kurganov and G. Petrova, {\em Commun. Math. Sci.}, 2007). The positivity of the
computed water height is ensured following (A.~Bollermann, S.~Noelle and
M.~Luk\'a\v{c}ov\'a, {\em Commun. Comput. Phys.}, 2010): The outgoing fluxes are
limited in case of draining cells.

\end{abstract}

\begin{keywords}
\smallskip
Hyperbolic systems of conservation and balance laws,
Saint-Venant system of shallow water equations,
finite volume methods,
well-balanced schemes,
positivity preserving schemes,
wet/dry fronts.
\smallskip

{\bf AMS subject classifications.} 76M12, 35L65
\end{keywords}

\section{Introduction}\label{intro}

We study numerical methods for the Saint-Venant system of shallow water equations \cite{SV}, which is widely used
for the flow of water in rivers or in the ocean. In one dimension, the Saint-Venant system reads:
\begin{eqnarray}
\left\{\begin{array}{l}\displaystyle{h_t+(hu)_x=0,}\\[0.3ex]
\displaystyle{(hu)_t +\Big(hu^2+\frac{1}{2}gh^2\Big)_x=-ghB_x},\end{array}\right.
\label{1d}
\end{eqnarray}
subject to the initial conditions
$$
h(x,0)=h_0(x),\quad u(x,0)=u_0(x),
$$
where $h(x,t)$ is the fluid depth, $u(x,t)$ is the velocity, $g$ is the gravitational constant, and the function
$B(x)$ represents the bottom topography, which is assumed to be independent of time $t$ and possibly discontinuous. The systems \eref{1d} is
considered in a certain spatial domain $X$ and if $X\neq\mathbb{R}$ the Saint-Venant system must be augmented
with proper boundary conditions.

In many applications, quasi steady solutions of the system \eref{1d} are to be captured using a (practically
affordable) coarse grid. In such a situation, small perturbations of steady states may be amplified by the scheme
and the so-called numerical storm can spontaneously develop \cite{NPPN}. To prevent it, one has to develop a well-balanced
scheme---a scheme that is capable of exactly balancing the flux and source terms so that ``lake at rest'' steady
states,
\begin{equation}
u=0, \quad w:=h+B={\rm Const}.
\label{sss}
\end{equation}
are preserved within the machine accuracy. Here, $w$ denotes the total water height or {\em free surface}.
Examples of such schemes can be found in
\cite{ABBKP,BNL,GPC2007,GHS,Jin,JW,KL,KP,LeV,LeVbook,NPPN,PS,RB,Rus1,Rus2,XS1,XS2}.

Another difficulty one often has to face in practice is related to the presence
of dry areas (island, shore) in the computational domain. As the eigenvalues of
the Jacobian of the fluxes in \eref{1d} are $u\pm\sqrt{gh}$, the system
\eref{1d} will not be strictly hyperbolic in the dry areas ($h=0$), and if due
to numerical oscillations $h$ becomes negative, the calculation will simply
break down. It is thus crucial for a good scheme to preserve the positivity of
$h$ (positivity preserving schemes can be found, e.g., in
\cite{ABBKP,GPC2007,BNL,KL,KP,PS,RB}).

We would also like to point out that when $h=0$ the ``lake at rest'' steady state \eref{sss} reduces to
\begin{equation}
hu=0, \quad h=0,
\label{dss}
\end{equation}
which can be viewed as a ``dry lake''. A good numerical scheme may be considered
``truly'' well-balanced when it is capable of exactly preserving both ``lake at
rest'' and ``dry lake'' steady states, as well as their combinations
corresponding to the situations, in which the domain $X$ is split into two
nonoverlapping parts $X_1$ (wet area) and $X_2$ (dry area) and the solution
satisfies \eref{sss} in $X_1$ and \eref{dss} in $X_2$.

\begin{figure}[htpb]
\centering
\includegraphics[width=5in,height = 2.5in]{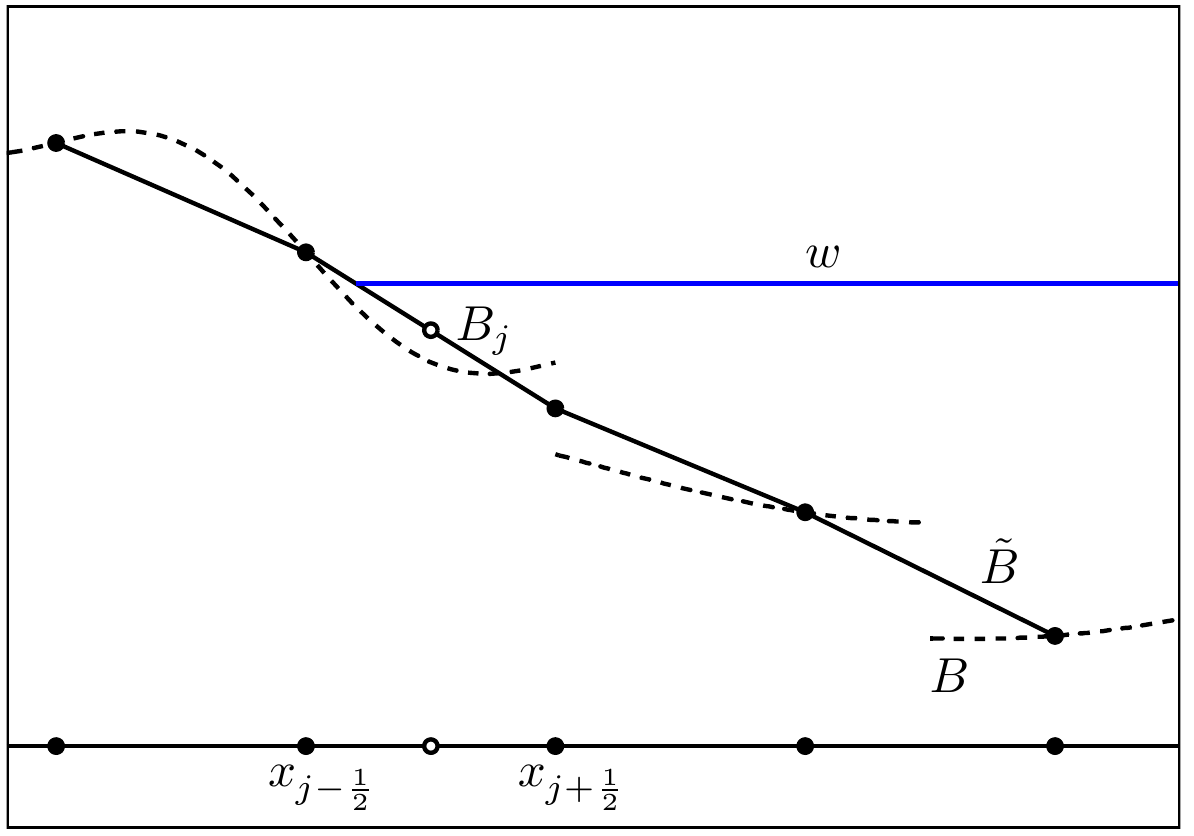}
%\caption{\sf ``Lake at rest'' steady state combined with dry boundaries.}
\caption{\sf ``Lake at rest'' steady state $w$ with dry boundaries upon a piecewise smooth topography $B$ (dashed line), which is reconstructed using piecewise linear, continuous $\tilde B$ (full line).}
\label{fig:wbrec}
\end{figure}

We focus on Godunov-type schemes, in which a numerical solution realized at a
certain time level by a global (in space) piecewise polynomial reconstruction,
is evolved to the next time level using the integral form of the system of
balance laws. In order to design a well-balanced scheme for \eref{1d}, it is
necessary that this reconstruction respects both the ``lake at rest'' \eref{sss}
and ``dry lake'' \eref{dss} steady-state solutions as well as their
combinations. On the other hand, to preserve positivity we have to make sure
that the reconstruction preserves a positive water height for all reconstructed
values. Both of this has been achieved by the hydrostatic reconstruction
introduced by Audusse et al. \cite{ABBKP}, based on a {\em discontinuous},
piecewise smooth discretisation of the bottom topography. In this paper, we
consider a {\em continuous}, piecewise linear reconstruction of the bottom. We propose
a piecewise linear reconstruction of the flow variables that also leads to a
well-balanced, positivity preserving scheme. The new reconstruction is based on
the proper discretization of a front cell in the situation like the one depicted
in Figure~\ref{fig:wbrec}. The picture depicts the real situation with a sloping
shore, and we see a discretization of the same situation that seems to be the
most suitable from a numerical perspective. We also demonstrate that the correct
handling of \eref{sss}, \eref{dss} and their combinations leads to a proper
treatment of non-steady states as well.

Provided the reconstruction preserves positivity, we can prove that the resulting central-upwind scheme is
positivity preserving. In fact, the proof from \cite{KP} carries over to the new scheme, but with a possibly
severe time step constraint. We therefore adopt a technique from \cite{BNL} and limit outgoing fluxes whenever the
so-called local draining time is smaller than the global time step. This approach ensures positive water heights
without a reduction of the global time step.

The paper is organized as follows. In \S\ref{sec:cu}, we briefly review the well-balanced positivity preserving
central-upwind scheme from \cite{KP}. A new positivity preserving reconstruction is presented in
\S\ref{sec:newrec}. The well-balancing and positivity preserving properties properties of the new scheme are
proven in \S\ref{sec:proofs}. Finally, we demonstrate the performance of the proposed method in
\S\ref{sec:exp}.

\section{A Central-Upwind Scheme for the Shallow Water Equations}\label{sec:cu}

Our work will be based on the central-upwind scheme proposed in \cite{KP}. We will therefore begin with a brief
overview of the original scheme.

We introduce a uniform grid $x_\alpha:=\alpha\dx$, with finite volume cells $I_j:=[x_\jmh,x_\jph]$ of length
$\dx$ and denote by $\overline\mU_j(t)$ the cell averages of the solution $\mU:=(w,hu)^T$ of \eref{1d} computed at
time $t$:
\begin{equation}
\overline\mU_j(t)\approx\frac{1}{\dx}\int\limits_{I_j}\mU(x,t)\,dx.
\label{eq:cellave}
\end{equation}
We then replace the bottom function $B$ with its continuous, piecewise linear approximation $\widetilde B$. To
this end, we first define
\begin{equation}
B_\jph:=\frac{B(x_\jph+0)+B(x_\jph-0)}{2},
\label{zv1}
\end{equation}
which in case of a continuous function $B$ reduces to $B_\jph=B(x_\jph)$, and then interpolate between these
points to obtain
\begin{equation}
\widetilde{B}(x)=B_\jmh+\left(B_\jph-B_\jmh\right)\cdot\frac{x-x_\jmh}{\dx},\quad x_\jmh\le x\le x_\jph.
\label{zv}
\end{equation}
From \eref{zv}, we obviously have
\begin{equation}
B_j:=\widetilde B(x_j)=\frac{1}{\dx}\int\limits_{I_j}\widetilde B(x)\,dx=\frac{B_\jph+B_\jmh}{2}.
\label{bot}
\end{equation}
The central-upwind semi-discretization of \eref{1d} can be written as the following system of time-dependent
ODEs:
\begin{equation}
\frac{d}{dt}\overline\mU_j(t)=-\frac{\mH_\jph(t)-\mH_\jmh(t)}{\dx}+\overline\mS_j(t),
\label{f1d}
\end{equation}
where $\mH_\jph$ are the central-upwind numerical fluxes and $\overline\mS_j$ is an appropriate discretization of
the cell averages of the source term:
\begin{equation}
\overline\mS_j(t)\approx\frac{1}{\dx}\int\limits_{I_j}\mS(\mU(x,t),B(x))\,dx,\quad\mS:=(0,-ghB_x)^T.
\label{dS}
\end{equation}
Using the definitions \eref{zv1} and \eref{bot}, we write the second component of the discretized source term
\eref{dS} as (see \cite{KL} and \cite{KP} for details)
\begin{equation}
 \overline\mS_j^{(2)}(t)
:=
 - g \overline h_j\frac{B_\jph-B_\jmh}{\dx}.
\label{s1}
\end{equation}

The central-upwind numerical fluxes $\mH_\jph$ are given by:
\begin{eqnarray}
\mH_\jph(t)&=&\frac{a^+_\jph\mF(\mU^-_\jph,B_\jph)-a^-_\jph\mF(\mU^+_\jph,B_\jph)}{a^+_\jph-a^-_\jph}\nonumber\\
&+&\frac{a^+_\jph a^-_\jph}{a^+_\jph-a^-_\jph}\left[\mU^+_\jph-\mU^-_\jph\right],
\label{nflux}
\end{eqnarray}
where we use the following flux notation:
\begin{equation}
\mF(\mU,B):=\left(hu,\frac{(hu)^2}{w-B}+\frac{g}{2}(w-B)^2\right)^T.
\label{eq:flux}
\end{equation}
The values $\mU^{\pm}_\jph =(w^{\pm}_\jph , h^{\pm}_\jph \cdot u^{\pm}_\jph ) $ represent the left and right values of the solution at point $x_\jph$ obtained by a
piecewise linear reconstruction
\begin{equation}
\widetilde{q}(x):=\overline q_j+(q_x)_j(x-x_j),\quad x_\jmh<x<x_\jph,
\label{r1}
\end{equation}
of $q$ standing for $w$ and $u$ respectively with $h^{\pm}_\jph = w^{\pm}_\jph - B_\jph$. To avoid the cancellation problem near dry areas,
we define the average velocity by
$$
\overline u_j:=\left\{\begin{array}{ll}(\overline {hu})_j/\overline h_j,&~\,\mbox{if}~\,\overline h_j\ge\epsilon,\\0,&~\mbox{otherwise}.
\end{array}\right.
$$
We choose $\epsilon=10^{-9}$ in all of our numerical experiments. This reconstruction will be second-order accurate if the approximate
values of the derivatives $(q_x)_j$ are at least first-order approximations of the corresponding exact derivatives. To ensure a
non-oscillatory nature of the reconstruction \eref{r1} and thus to avoid spurious oscillations in the numerical solution, one has to
evaluate $(q_x)_j$ using a nonlinear limiter. From the large selection of the limiters readily available in the literature (see, e.g.,
\cite{GR2,Kr,vLeV,LeVbook,LieNoe2,NT,Swe}), we chose the generalized minmod limiter (\cite{vLeV,LieNoe2,NT,Swe}):
\begin{equation}
(q_x)_j={\rm minmod}\left(\theta\frac{\overline q_j-\overline q_{j-1}}{\dx},\,\frac{\overline q_{j+1}-\overline q_{j-1}}{2\dx},\,
\theta\frac{\overline q_{j+1}-\overline q_j}{\dx}\right),\quad\theta\in[1,2],
\label{minmod}
\end{equation}
where the minmod function, defined as
\begin{equation}
{\rm minmod}(z_1,z_2,...):=\left\{\begin{array}{lc}\!\!\!\min_j\{z_j\}, & ~~\mbox{if} ~~z_j>0 ~~\forall j,\\
\!\!\!\max_j\{z_j\}, & ~~\mbox{if} ~~z_j<0 ~~\forall j,\\ \!\!\!0, & ~\mbox{otherwise},\end{array}\right.
\label{mm}
\end{equation}
is applied in a componentwise manner, and $\theta$ is a parameter affecting the numerical viscosity of the
scheme. It is shown in \cite{KP} that this procedure (as well as any alternative ``conventional'' reconstruction,
including the simplest first-order piecewise constant one, for which $(w_x)_j\equiv\mo$) might produce negative
values $h_{\jph}^\pm$ near the dry areas (see \cite{KP}). Therefore, the reconstruction \eref{r1}--\eref{mm} must
be corrected there. The correction algorithm used in \cite{KP} restores positivity of the reconstruction depicted in Figure~\ref{fig:kp}, but
destroys the well-balancing property. This is explained in \S\ref{sec:newrec}, where we propose an alternative
positivity preserving reconstruction, which is capable of exactly preserving the ``lake at rest'' and the ``dry
lake'' steady states as well as their combinations.

Finally, the local speeds $a^\pm_\jph$ in \eref{nflux} are obtained using the eigenvalues of the Jacobian
$\frac{\partial\mF}{\partial\mU}$ as follows:
\begin{align}
 a_\jph^+
& =
 \max\left\{u_\jph^++\sqrt{gh_\jph^+}\,,\,u_\jph^-+\sqrt{gh_\jph^-}\,,\,0\right\},
\label{lsp}
\\
 a_\jph^-
& =
 \min\left\{u_\jph^+-\sqrt{gh_\jph^+}\,,\,u_\jph^--\sqrt{gh_\jph^-}\,,\,0\right\}.
\label{lsp1}
\end{align}
Note that for $\overline\mU_j$, $\mU^\pm_\jph$ and $a^\pm_\jph$, we dropped the dependence of $t$ for simplicity.

As in \cite{KP}, in our numerical experiments, we use the third-order strong
stability preserving Runge-Kutta (SSP-RK) ODE solver (see \cite{GST} for
details) to numerically integrate the ODE system \eref{f1d}. The timestep is
restricted by the standard CFL condition,
\begin{equation}
\label{CFL}
\CFL := \frac\dt\dx \max\limits_j |a^\pm_{j+\frac12}| \;\; \le \;\; \frac12
\end{equation}
For the examples of the present paper, results of the second and third order
SSP-RK solvers are almost undistinguishable.

\section{A New Reconstruction at the Almost Dry Cells}\label{sec:newrec}

\begin{figure}[htpb]
\centering
\includegraphics[width=5in,height = 2.5in]{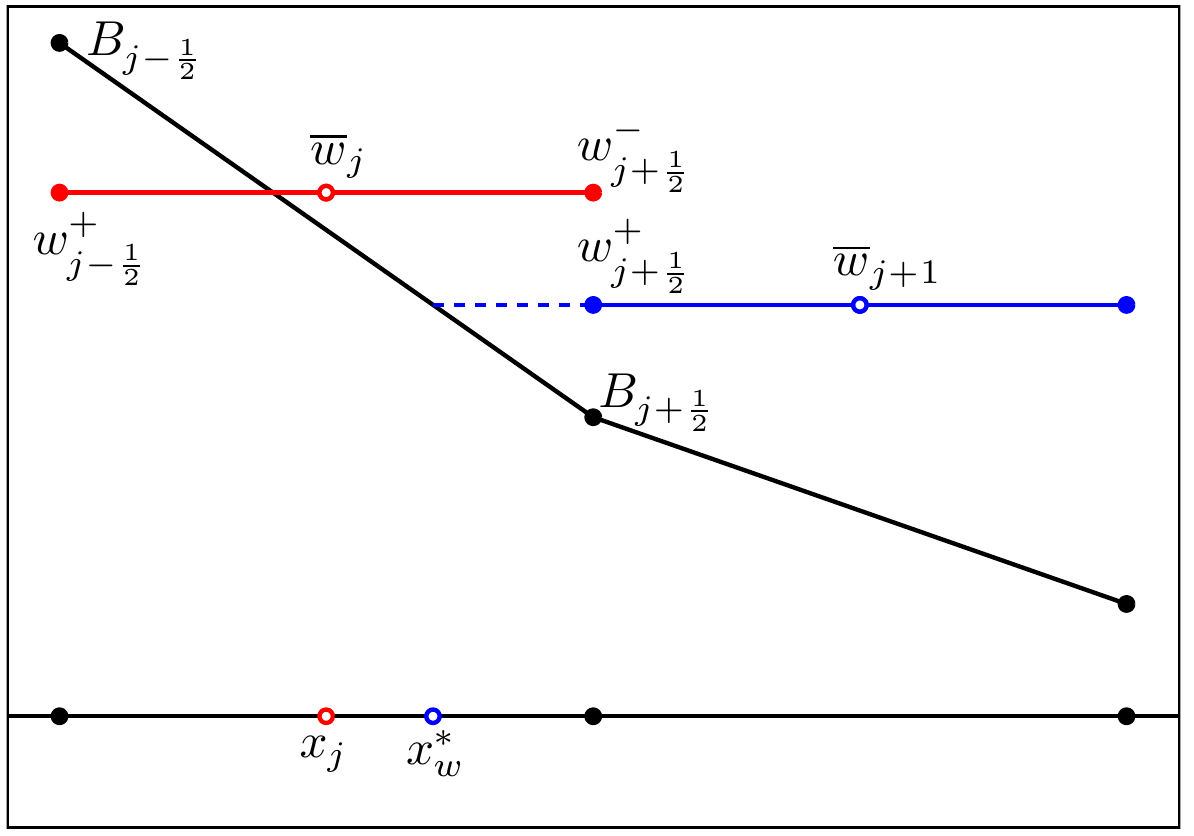}
\caption{Wrong approximations of the wet/dry front by the piecewise constant reconstruction.}
\label{fig:pconstrec}
\end{figure}

In the presence of dry areas, the central-upwind scheme described in the previous section may create negative
water depth values at the reconstruction stage. To understand this, one may look at Figure \ref{fig:pconstrec}, where we
illustrate the following situation: The solution satisfies \eref{sss} for $x>\xs_w$ (where $\xs_w$ marks the
waterline) and \eref{dss} for $x<\xs_w$. Notice that cell $j$ is a typical almost dry cell and the use of the
(first-order) piecewise constant reconstruction clearly leads to appearance of negative water depth values there.
Indeed, in this cell the total amount of water is positive and therefore $\overline w_j>B_j$, but clearly
$\overline w_j<B_{j-\frac{1}{2}}$ and thus $h_{j-\frac{1}{2}}<0$.

\begin{figure}[htpb]
\centering
\includegraphics[width=5in,height = 2.5in]{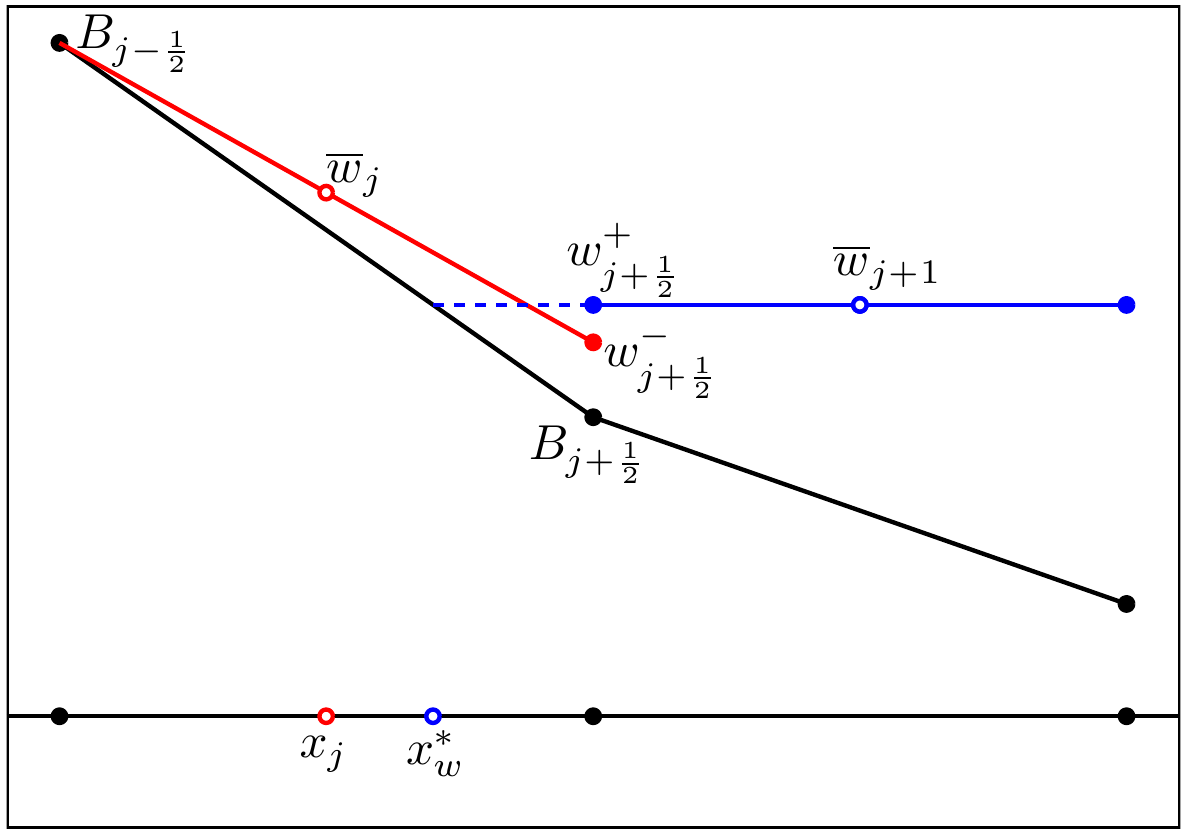}
\caption{Approximations of the wet/dry front by the positivity preserving but unbalanced piecewise linear
reconstruction from \hbox{ \cite{KP} }.}
\label{fig:kp}
\end{figure}

It is clear that replacement of the first-order piecewise constant reconstruction with a conventional
second-order piecewise linear one will not guarantee positivity of the computed point values of $h$. Therefore,
the reconstruction in cell $j$ may need to be corrected. The correction proposed in \cite{KP} will solve the
positivity problem by raising the water level at one of the cell edges to the level of the bottom function there
and lowering the water level at the other edge by the same value (this procedure would thus preserve the amount
of water in cell $j$). The resulting linear piece is shown in Figure \ref{fig:kp}. Unfortunately, as one may
clearly see in the same figure, the obtained reconstruction is not well-balanced since the reconstructed values
$w_\jph^-$ and $w_\jph^+$ are not the same.

\begin{figure}[htpb]
\centering
\includegraphics[width=5in,height = 2.5in]{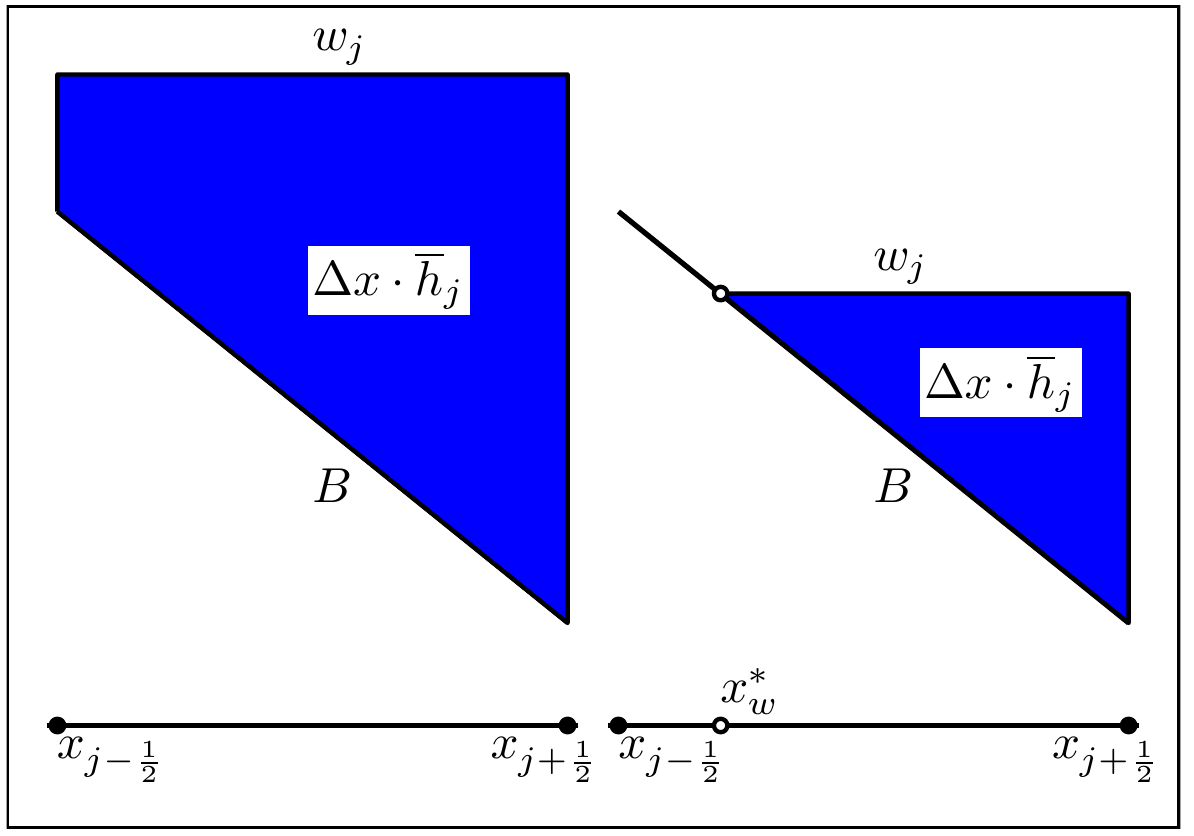}
\caption{\sf Computation of $w_j$. Left: Fully flooded cell; Right: Partially flooded cell.}
\label{fig:compute_w}
\end{figure}

Here, we propose an alternative correction procedure, which will be both positivity
preserving and well-balanced even in the presence of dry areas.
This correction bears some similarity to the reconstruction near dry fronts of
depth-averaged granular avalanche
models in \cite{TNGH}. However, in \cite{TNGH} the authors tracked a front
running down the terrain, and did not treat well-balancing of equilibrium states.
Let us assume that at a certain time level all computed values
$\overline w_j\ge B_j$ and the slopes $(w_x)_j$ and $(u_x)_j$ in the piecewise linear
reconstruction \eref{r1} have been
computed using some nonlinear limiter as it was discussed in \S\ref{sec:cu} above.
We also assume that at some almost dry cell $j$,
\begin{equation}
%B_\jph<\overline w_j<B_\jmh
B_\jmh>\overline w_j>B_\jph
\label{3.1}
\end{equation}
(the case $B_\jmh<\overline w_j<B_\jph$ can obviously be treated in a symmetric way) and that the reconstructed values of $w$ in cell $j+1$
satisfy
\begin{equation}\label{eqn:dry-cond}
w_\jph^+>B_\jph \quad \hbox{and} \quad w_{j+\frac{3}{2}}^->B_{j+\frac{3}{2}},
\end{equation}
that is, cell $j+1$ is fully flooded. This means that cell $j$ is located
near the dry boundary (mounting shore), and we design a well-balanced reconstruction correction procedure for cell $j$ in the following way:

We begin by computing the free surface in cell $j$ (denoted by $w_j$), which represents the average total water
level in ({\em the flooded parts} of) this cell assuming that the water is at rest.
The meaning of this formulation becomes clear from Figure~\ref{fig:compute_w}. We always
choose $w_j$ such that the area enclosed between the line with height $w_j$ and the
bottom line equals the amount of water given by $\dx\cdot\overline h_j$, where
$\overline h_j:=\overline w_j-B_j$. The resulting area is either a trapezoid
(if cell $j$ is a fully flooded cell as in Figure~\ref{fig:compute_w} on the left) or
a triangle (if cell $j$ is a partially flooded cell as in Figure~\ref{fig:compute_w} on the right),
depending on $\overline h_j$ and the bottom slope $(B_x)_j$.

So if the cell $j$ is a fully flooded cell, i.e. $~\overline h_j\ge\frac{\dx}{2}\,\left|(B_x)_j\right|$, the free surface $w_j(x)$ is defined as
\[ w_j(x) = \overline w_j, \]
otherwise the free surface is a continuous piecewise linear function given by
\begin{equation}
w_j(x)=\begin{cases}B_j(x),&\text{if $x<\xs_w$},\\
w_j,&\text{otherwise},\end{cases}
\label{eq:wj}
\end{equation}
where $\xs_w$ is the boundary point separating the dry and wet parts in the cell $j$. It can be determined by the mass conservation,
\begin{eqnarray*}
\dx\cdot\overline h_j&=&\int_{x_\jmh}^{x_\jph}(w_j(x)-B_j(x))dx = \int_{\xs_w}^{x_\jph}(w_j-B_j(x))dx \\
&=&\frac{\dxs_w}{2}\,(w_j-B_\jph)
=\frac{\dxs_w}{2}\,(B(\xs_w)-B_\jph)=-\frac{(\dxs_w)^2}{2}(B_x)_j,
\end{eqnarray*}
where $\dxs_w=x_\jph-\xs_w$, thus
\begin{equation}
\dxs_w=\sqrt{\frac{2\dx\overline h_j}{-(B_x)_j}}=
\sqrt{\frac{2\overline h_j}{B_\jmh-B_\jph}}~\dx,
\label{eq:dxw}
\end{equation}
resulting in the free surface $w_j$ formula for the wet/dry cells,
 \begin{equation}
w_j= B_{\jph} + \sqrt{{2\overline h_j}|B_\jmh-B_\jph|}
\label{eq:wjdry}
\end{equation}
Note that the limit for the distinction of cases in \eref{eq:wj} is determined from the area of the triangle
between the bottom line and the horizontal line at the level of $B_\jmh$. We also note that if cell $j$
satisfies \eref{3.1}, then it is clearly a partially flooded cell (like the one shown in Figure~\ref{fig:compute_w} on
the right) with $\dxs_w<\dx$.

\smallskip
\begin{rem}
We would like to emphasize that if cell $j$ is fully flooded, then the free surface is represented by the cell
average $\overline w_j$ (see the first case in equation {\rm\eref{eq:wj}}), while if the cell is only partially
flooded, $\overline w_j$ does not represent the free surface at all (see, e.g., Figure {\rm\ref{fig:pconstrec}}). Thus,
in the latter case we need to represent the free surface with the help of another variable $w_j\neq\overline w_j$
(see the second case in {\rm\eref{eq:wj}}), which is only defined on the wet part of cell $j$, $[\xs_w,x_\jph]$,
and thus stays above the bottom function $B$, see Figure {\rm\ref{fig:compute_w}} (right).
\end{rem}

\begin{figure}[htpb]
\centering
\includegraphics[width=5in,height = 2.5in]{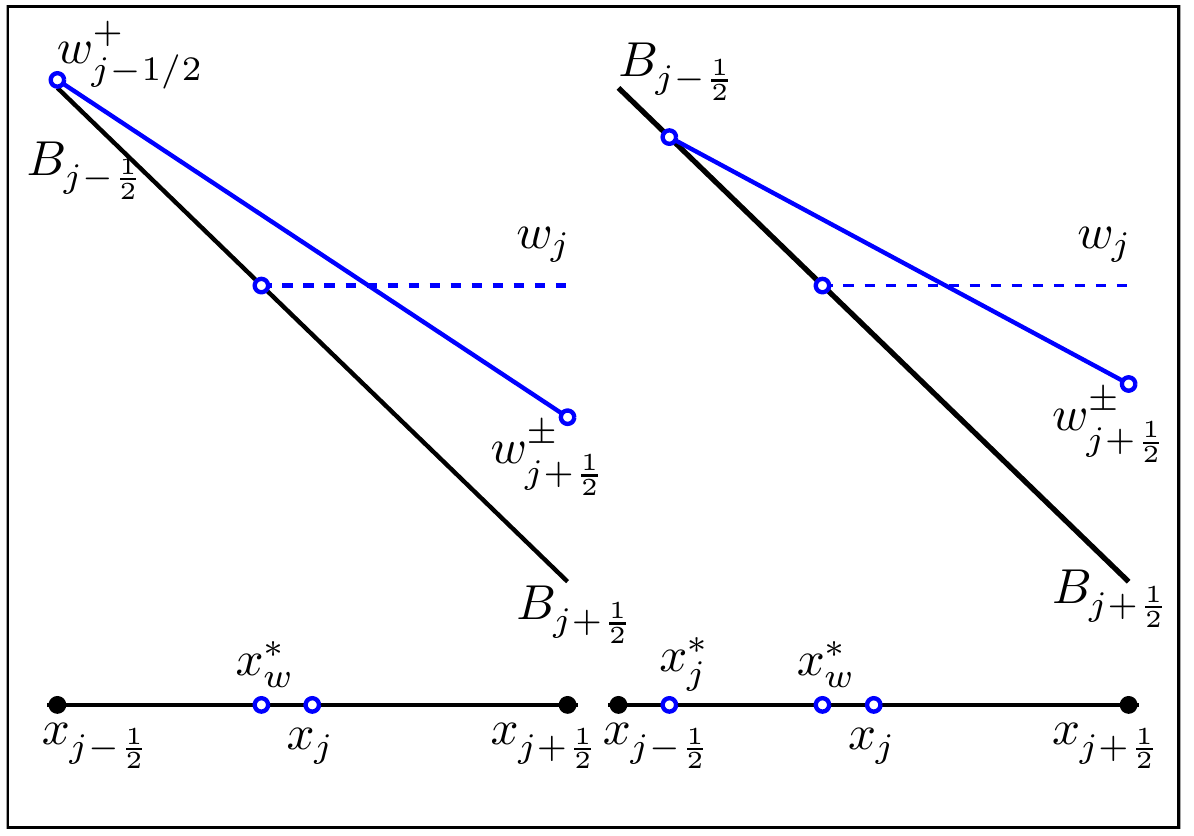}
\caption{\sf Conservative reconstruction of $w$ at the boundary with the fixed value $w^+_\jph$. Left: Linear
reconstruction with nonnegative $h^+_\jmh$; Right: Two linear pieces with $h^+_\jmh=0$.}
\label{fig:LinearRec}
\end{figure}

\smallskip
We now modify the reconstruction of $h$ in the partially flooded cell $j$ to ensure the well-balanced property. To this
end, we first set $w_\jph^-=w_\jph^+$ (which immediately implies that
$h_\jph^-:=w_\jph^--B_\jph=w_\jph^+-B_\jph=:h_\jph^+$) and determine the reconstruction of $w$ in cell $j$ via
the conservation of $\overline h_j$ in this cell. We distinguish between the following two possible cases. If the amount of water in cell $j$ is sufficiently large (as in the case illustrated in Figure \ref{fig:LinearRec} on the left), there is a unique $h_\jmh^+ \geq 0$ satisfying
\begin{equation}
\overline h_j=\frac{1}{2}(h^-_\jph+h^+_\jmh).
\label{3.4}
\end{equation}
From this we obtain $w_\jmh^+=h^+_\jmh+B_\jmh$, and thus the well-balanced reconstruction in cell $j$ is completed.

If the value of $h_\jmh^+$, computed from the conservation requirement \eref{3.4} is negative, we replace a
linear piece of $w$ in cell $j$ with two linear pieces as shown in Figure \ref{fig:LinearRec} on the right. The
breaking point between the ``wet'' and ``dry'' pieces will be denoted by $\xs_j$ and it will be determined from
the conservation requirement, which in this case reads
\begin{equation}
\dx\cdot\overline h_j=\frac{\dxs_j}{2}\,h^-_\jph,
\label{eq:reccons}
\end{equation}
where
$$
\dxs_j=\big|x_\jph-\xs_j\big|.
$$

Combining the above two cases, we obtain the reconstructed value
\begin{equation}
h^+_\jmh=\max\left\{0,\,2\overline h_j-h^-_\jph\right\}.
\label{eq:rech}
\end{equation}
We also generalize the definition of $\dxs_j$ and set
\begin{equation}
\dxs_j:=\dx\cdot\min\left\{\frac{2\,\overline h_j}{h^-_\jph},\,1\right\},
\label{eq:dxr}
\end{equation}
which will be used in the proofs of the positivity and well-balancing of the resulting central-upwind scheme in
\S\ref{sec:proofs}.
We summarize the wet/dry reconstruction in the following definition:
\begin{defn}\label{def:rec}(wet/dry reconstruction)
For the sake of clarity, we denote the left and right values of the piecewise linear reconstruction \eqref{r1} -- \eqref{mm} by $\tilde{\mU}^{\pm}_\jph =(\tilde w^{\pm}_\jph , \tilde h^{\pm}_\jph \cdot \tilde u^{\pm}_\jph )$. The purpose of this definition is to define the final values
 $\mU^{\pm}_\jph =(w^{\pm}_\jph , h^{\pm}_\jph \cdot u^{\pm}_\jph )$, which are modified by the wet/dry
 reconstruction.
%\\[0.5ex]
%{\bf Case 1:} $ \bar w_j \geq \max(B_{j\pm\frac12})$
%\\
%{\bf Case 1a:} $\tilde w_{j\pm\frac12} \geq \max(B_{j\pm\frac12})$
%\\
%The cell is fully flooded, and we set $\mU^{\pm}_\jph := \tilde{\mU}^{\pm}_\jph$.
%\\
%{\bf Case 1b:} $\tilde w_{j-\frac12} < \max(B_{j-\frac12})$ or
%$\tilde w_{j+\frac12} < \max(B_{j+\frac12})$
%\\
%There is enough water to flood the cell, and as in \cite{KP} we redistribute the water via
%\begin{eqnarray*}
%&&\mbox{If}~~w^-_\jph<B_\jph,~~\mbox{then set}~~(w_x)_j:=\frac{B_\jph-\overline w_j}{\dx/2},\nonumber\\
%&&\Longrightarrow~~w^-_\jph=B_\jph,~~ w^+_\jmh=2\overline w_j-B_\jph;
%\end{eqnarray*}
%and
%\begin{eqnarray*}
%&&\mbox{If}~~w^+_\jmh<B_\jmh,~~\mbox{then set}~~(w_x)_j:=\frac{\overline w_j-B_\jmh}{\dx/2},\nonumber\\
%&&\Longrightarrow~~w^-_\jph=2\overline w_j-B_\jmh, ~~ w^+_\jmh=B_\jmh.
%\end{eqnarray*}
%\\
%{\bf Case 2:} $ \min(B_{j\pm\frac12}) < \bar w_j < \max(B_{j\pm\frac12})$
%\\
%The cell is partially flooded.
%\\
%{\bf Case 2a:}  $\tilde w_{j-\frac12} < B_{j-\frac12}$ and $\tilde w_{j+\frac12} \geq B_{j-\frac12}$
%\\
%{\bf Case 2a1:} cell $(j+1)$ is fully flooded.
%\\
%Define $h_\jph^-$  using \eref{eq:reccons} and $\dxs_j$ by \eref{eq:dxr}.
%\\
%{\bf Case 2a2:}
%otherwise set $h_\jph^-:=w_j-B_\jph$~\eref{eq:wjdry}
%and $\dxs_j:=\dxs_w$~\eref{eq:dxw}. Note that this situation is not
%generic and may occur only in the under-resolved computations.
%\\
%{\bf Case 2b:}  $\tilde w_{j-\frac12} \geq B_{j-\frac12}$ and $\tilde w_{j+\frac12} < B_{j-\frac12}$
%\\
%Analogous to Case 2a.

\begin{description}
  \item[Case 1.]  $ \bar w_j \geq B_{j-\frac12}$ and  $ \bar w_j \geq B_{j+\frac12}$: there is enough water to flood the cell for flat lake.
  \begin{description}
    %\item[1A.] $\tilde w_{j\pm\frac12} \geq \max(B_{j\pm\frac12})$:
    \item[1A.] $\tilde w_{j-\frac12}^+ \geq B_{j-\frac12}$ and $\tilde w_{j+\frac12}^- \geq B_{j+\frac12}$  :
       the cell is fully flooded, and we set $\mU^{\pm}_\jph := \tilde{\mU}^{\pm}_\jph$.
    \item[1B.]otherwise, as in \cite{KP} we redistribute the water via
         \begin{eqnarray*}
           &&\mbox{If}~~\tilde w^-_\jph<B_\jph,~~\mbox{then set}~~(w_x)_j:=\frac{B_\jph-\overline w_j}{\dx/2},\nonumber\\
           &&\Longrightarrow~~w^-_\jph=B_\jph,~~ w^+_\jmh=2\overline w_j-B_\jph;
         \end{eqnarray*}
         and
         \begin{eqnarray*}
           &&\mbox{If}~~\tilde w^+_\jmh<B_\jmh,~~\mbox{then set}~~(w_x)_j:=\frac{\overline w_j-B_\jmh}{\dx/2},\nonumber\\
           &&\Longrightarrow~~w^-_\jph=2\overline w_j-B_\jmh, ~~ w^+_\jmh=B_\jmh.
         \end{eqnarray*}
   \end{description}
     \item[Case 2.]  $B_\jmh>\overline w_j>B_\jph$: the cell is possible partially flooded.
             \begin{description}
               \item[2A.] $\tilde w_\jph^+>B_\jph$ and $\tilde w_{j+\frac{3}{2}}^->B_{j+\frac{3}{2}}$, i.e., cell $j+1$ is fully flooded and
                   $w_\jph^+ = \tilde w_\jph^+$ . Define $w_\jph^- = w_\jph^+$ and $h_\jph^- = w_\jph^- - B_\jph$ .
                   \begin{description}
                   \item[2A1.] $2 \overline h_j - h_\jph^- \geq 0$, the amount of water in cell $j$ is sufficiently large, we set  $ h^+_\jmh = 2 \overline h_j - h_\jph^-$, so $w_\jmh^+=h^+_\jmh+B_\jmh$
                   \item[2A2.] otherwise set $h^+_\jmh = 0, w_\jmh^+=B_\jmh$ and $\dxs_j$ as in \eref{eq:dxr}.
                   \end{description}
               \item[2B.]otherwise set $h_\jph^-:=w_j-B_\jph$~\eref{eq:wjdry} and $\dxs_j:=\dxs_w$~\eref{eq:dxw}. Note that this situation is not generic and may occur only in the under-resolved computations.
             \end{description}
    \item[Case 3.]  $B_\jmh<\overline w_j<B_\jph$: analogous to {\rm \textbf{Case 2}}.
\end{description}
\end{defn}

%%%%%%%%%%%%%%%%%%%%%%%%%%%%%%%%%%%%%%%%%%%%%%%%%%%%%%%%%%%%%%%%%

\section{Positivity Preserving and Well-Balancing}\label{sec:proofs}

In the previous section, we proposed a new spatial reconstruction for wet/dry
cell. In this section, we will implement a time-quadrature for the fluxes at
wet/dry boundaries developed in \cite{BNL}. It cuts off the space-time flux
integrals for partially flooded interfaces. Then we prove that the resulting
central-upwind scheme is positivity preserving and well-balanced under the
standard CFL condition \eqref{CFL}.

We begin by studying the positivity using a standard time integration of the
fluxes. The following lemma shows that for explicit Euler time stepping,
positivity cannot be guaranteed directly under a CFL condition such as
\eqref{CFL}.

\smallskip

\begin{lem}\label{LEMpos}
{\rm\eref{f1d}--\eref{lsp1}} with the piecewise linear reconstruction
{\rm\eref{r1}} corrected according to the procedure described in
\S\ref{sec:newrec}. Assume that the system of ODEs {\rm\eref{f1d}} is solved by
the forward Euler method and that for all $j$, $\overline h^n_j\ge 0$. Then
\\[1ex]
(i) $\overline h^{n+1}_j\ge 0$ for all $j$ provided that
\begin{equation}
 \dt
\le
 \min_j\left\{\frac{\dxs_j}{2a_j}\right\}, \qquad
 a_j
:=
 \max\{a_\jph^+,-a_\jph^-\}.
\label{CFLstar}
\end{equation}
(ii) Condition \eqref{CFLstar} cannot be guaranteed by any finite positive CFL
condition \eqref{CFL}.
\end{lem}

\smallskip

\noindent {\bf Proof:} (i) For the fully flooded cells with $\dxs_j=\dx$, the proof of Theorem 2.1 in \cite{KP} still holds.
Therefore, we will only consider partially flooded cells like the one shown in Figure \ref{fig:LinearRec}. First, from
\eref{eq:reccons} we have that in such a cell $j$ the cell average of the water depth at time level $t=t^n$ is
\begin{equation}
\overline h_j^n=\frac{\dxs_j}{2\dx}\,h^-_\jph,
\label{4.1}
\end{equation}
and it is evolved to the next time level by applying the forward Euler temporal discretization to the first
component of \eref{f1d}, which after the subtraction of the value $B_j$ from both sides can be written as
\begin{equation}
\overline h_j^{n+1}=\overline h_j^n-\lambda\left(\mH^{(1)}_\jph-\mH^{(1)}_\jmh\right),\quad
\lambda:=\frac{\dt}{\dx},
\label{FE}
\end{equation}
where the numerical fluxes are evaluated at time level $t=t^n$. Using \eref{nflux} and the fact that by
construction $w_\jph^+-w_\jph^-=h_\jph^+-h_\jph^-$, we obtain:
\begin{equation}
\mH^{(1)}_\jph=\frac{a_\jph^+(hu)_\jph^--a_\jph^-(hu)_\jph^+}{a_\jph^+-a_\jph^-}+
\frac{a_\jph^+a_\jph^-}{a_\jph^+-a_\jph^-}\Big[h_\jph^+-h_\jph^-\Big].
\label{4.3}
\end{equation}
Substituting \eref{4.1} and \eref{4.3} into \eref{FE} and taking into account the fact that in this cell
$h_\jmh^+=0$, we arrive at:
\begin{eqnarray}
\overline h_j^{n+1}&=&
\left[\frac{\dxs_j}{2\dx}-\lambda a_\jph^+\left(\frac{u_\jph^--a_\jph^-}{a_\jph^+-a_\jph^-}\right)\right]h_\jph^-
\nonumber\\
&-&\lambda a_\jph^-\left(\frac{a_\jph^+-u_\jph^+}{a_\jph^+-a_\jph^-}\right)h_\jph^++
\lambda a_\jmh^+\left(\frac{u_\jmh^--a_\jmh^-}{a_\jmh^+-a_\jmh^-}\right)h_\jmh^-,
\label{eq:hupdate}
\end{eqnarray}
Next, we argue as in \cite[Theorem 2.1]{KP} and show that $\overline h_j^{n+1}$
is a linear combination of the three values, $h^\pm_\jph$ and $h^-_\jmh$ (which
are guaranteed to be nonnegative by our special reconstruction procedure) with
nonnegative coefficients. To this end, we note that it follows from \eref{lsp}
and \eref{lsp1} that $a_\jph^+\ge 0$, $a_\jph^-\le 0$, $a_\jph^+-u_\jph^+\ge 0$,
and $u_\jph^--a_\jph^-\ge 0$, and hence the last two terms in \eref{eq:hupdate}
are nonnegative. By the same argument,
$0\le\frac{a_\jmh^+-u_\jmh^+}{a_\jmh^+-a_\jmh^-}\le 1$ and
$0\le\frac{u_\jph^--a_\jph^-}{a_\jph^+-a_\jph^-}\le 1$, and thus the first term
in \eref{eq:hupdate} will be also nonnegative, provided the CFL restriction
\eref{CFLstar} is satisfied. Therefore, $\overline h_j^{n+1}\ge 0$, and part (i)
is proved.

In order to show part (ii) of the lemma, we compare the CFL-like conndition
\eqref{CFLstar} with the standard CFL condition \eqref{CFL},
\begin{align}
 \CFL^*
&:=
 \dt \; \max\limits_j \left( \frac{|a_j|}{\dx_j^*} \right) \;
% =
 % \frac\dt\dx \; \max\limits_j \left( |a_j| \frac{\dx}{\dx_j^*} \right) \;
=
 \max\limits_j \left(  \frac{|a_j|}{\max\limits_i |a_i|} \frac{\dx}{\dx_j^*} \right) \; \CFL
\end{align}
We note that depending on the water level $w_j$ in the partially flooded cell,
$\Delta x_j^*$ can be arbitrarily small, so there is no upper bound of $\CFL^*$
in terms of $\CFL$.
$\hfill\Box$

\smallskip

Part (ii) of Lemma~\ref{LEMpos} reveals that one might obtain a
serious restriction of the timestep in the presence of partially flooded cells.
We will now show how to overcome this restriction using the draining time
technique developed in \cite{BNL}.

We start from the equation \eqref{FE} for the water height and look for a suitable modification
 of the update such that the water height remains positive,
\begin{equation*}
\overline h_j^{n+1}=\overline h_j^n-\dt\,\frac{\mH^{(1)}_\jph-\mH^{(1)}_\jmh}{\dx}\ge 0.
\end{equation*}
As in \cite{BNL}, we introduce the {\em draining time step}
\begin{equation}
\label{eq:draining-time}
\dtd:=\frac{\dx\overline h_j^n}{\max(0,\mH^{(1)}_\jph)+\max(0,-\mH^{(1)}_\jmh)},
\end{equation}
which describes the time when the water contained in cell $j$ in the beginning of the time step has left via the
outflow fluxes. We now replace the evolution step \eref{FE} with
\begin{equation}
\label{eq:fvh}
\overline h_j^{n+1}=\overline h_j^n-\frac{\dt_\jph\mH^{(1)}_\jph-\dt_\jmh\mH^{(1)}_\jmh}{\dx},
\end{equation}
where we set the effective time step on the cell interface as
\begin{equation}
\label{eq:effect-time}
\dt_\jph=\min(\dt,\dt_i^{\text{drain}}),
\quad
i=j+\frac12-\frac{\sgn\left(\mH^{(1)}_\jph\right)}{2}.
\end{equation}
The definition of $i$ selects the cell in upwind direction of the edge. We would
like to point out that the modification of flux is only active in cells which
are at risk of running empty during the next time step. It corresponds to the
simple fact that there is no flux out of a cell once the cell is empty. The
positivity based on the draining time is summarized as the following
theorem, which we proved in \cite{BNL}. Note that in contrast to Lemma
\ref{LEMpos}, the timestep is now uniform under the CFL condition \eqref{CFL}:

\begin{theorem}\label{THMpos}
Consider the update \eqref{eq:fvh} of the water height with
fluxes with the help of the draining time
\eqref{eq:draining-time}. Assume that the initial height
$\overline h^n_j$ is non-negative for all $j$. Then the height
remains nonnegative,
\begin{equation}
 \overline h^{n+1}_j \geq 0 \quad \mathrm{for \; all} \;\; j.
\end{equation}
provided that the standard CFL condition \eqref{CFL} is satisfied.
\end{theorem}

To guarantee well-balancing, we have to make sure that the gravity driven part of the momentum flux
$\mH^{(2)}_\jph$ cancels the source term $\mS^{(2)}_\jph$, in a lake at rest situation.
To this end, we follow \cite{BNL} and split the momentum flux $\mF^{(2)}(\mU)$ in its advective and
gravity driven parts:
$$
 \mF^{(2),\rm a}(\mU):=\frac{(hu)^2}{w-B}~~\text{and}~~~
 \mF^{(2),\rm g}(\mU):=\frac{g}{2}(w-B)^2,
$$
respectively. For convenience, we will denote $w-B$ by $h$ in the following.
The corresponding advective and gravity driven parts of the central-upwind fluxes then read
$$
\mH^{(2),\rm g}_\jph(t)=\frac{a^+_\jph\mF^{(2),\rm g}(\mU^-_\jph)-a^-_\jph\mF^{(2),\rm g}(\mU^+_\jph)}{a^+_\jph-a^-_\jph}
+\frac{a^+_\jph a^-_\jph}{a^+_\jph-a^-_\jph}\left[\mU^{(2),+}_\jph-\mU^{(2),-}_\jph\right],
$$
and
$$
\mH^{(2),\rm a}_\jph(t)=\frac{a^+_\jph\mF^{(2),\rm a}(\mU^-_\jph)-a^-_\jph\mF^{(2),\rm a}(\mU^+_\jph)}{a^+_\jph-a^-_\jph},
$$
The above fluxes adds up to the following modified update of the momentum:
\begin{equation}
 (\overline{hu})_j^{n+1}=(\overline{hu})_j^{n}-\frac{\dt_\jph\mH^{(2),\rm a}_\jph-\dt_\jmh\mH^{(2),\rm a}_\jmh}{\dx}
-\dt\left(\frac{\mH^{(2),\rm g}_\jph-\mH^{(2),\rm g}_\jmh}{\dx}+\overline\mS_j^{(2),n}\right).
\label{eq:fv-m}
\end{equation}

This modified finite volume scheme  {\rm\eref{eq:fvh}} and {\rm\eref{eq:fv-m}}  ensures the well-balancing property
even in the presence of dry areas, as we will
show in Theorem~\ref{THMwb}.

\smallskip
\begin{theorem}\label{THMwb}
Consider the system {\rm\eref{1d}} and the fully discrete central-upwind scheme {\rm\eref{eq:fvh}} and {\rm\eref{eq:fv-m}}. Assume
that the numerical solution $\mU(t^n)$ corresponds to the steady state which is a combination of the ``lake at
rest'' {\rm\eref{sss}} and ``dry lake'' {\rm\eref{dss}} states in the sense that for all $w_j$ defined in
{\rm\ref{eq:wj}}, $w_j={\rm Const}$ and $u=0$ whenever $h_j>0$. Then $\mU(t^{n+1})=\mU(t^n)$, that is, the scheme
is {\em well-balanced}.
\end{theorem}

\smallskip
{\bf Proof:} We have to show that in all cells the fluxes and the source term
discretization cancel exactly. First, we mention the fact that the
reconstruction procedure derived in \S\ref{sec:newrec} preserves both the ``lake
at rest'' and ``dry lake'' steady states and their combinations. For all cells
where the original reconstruction is not corrected, the resulting slopes are
obviously zero and therefore $w^\mp_\jpmh=w_j$ there. As $hu=0$ in all cells,
the reconstruction for $hu$ obviously reproduces the constant point values
$(hu)^\mp_\jpmh=0,~\forall j$, resulting that the draining time is equal to the
global time step, i.e., $\dtd = \Delta t$.

We first analyze the update of the free surface using \eref{eq:fvh}. The first component of flux \eref{nflux}
is
$$
\mH^{(1)}_\jph=\frac{a_\jph^+(hu)_\jph^--a_\jph^-(hu)_\jph^+}{a_\jph^+-a_\jph^-}+
\frac{a_\jph^+a_\jph^-}{a_\jph^+-a_\jph^-}\Big[(h+B)_\jph^+-(h+B)_\jph^-\Big]=0,
$$
as $B_\jph^+=B_\jph^-$, $h_\jph^+=h_\jph^-$  and $(hu)_\jph^+=(hu)_\jph^-=0$. This gives
$$
\overline{w}_j^{n+1} = \overline{h}_j^{n+1} + B_j = \overline{h}_j^{n}+B_j = \overline{w}_j^n
$$

Secondly, we analyze the update of the momentum using \eref{eq:fv-m}.  Using the same argument and setting $u^\pm_\jph=0$ at the
points $x=x_\jph$ where $h^+_\jph=h^-_\jph=0$, for the second component we obtain
\begin{eqnarray*}
\mH^{{(2),\rm a}}_\jph+\mH^{{(2),\rm g}}_\jph&=&\frac{a_\jph^+\left(hu^2\right)_\jph^-
-a_\jph^-\left(hu^2\right)_\jph^+}{a_\jph^+-a_\jph^-}+\frac{a_\jph^+\left(\frac{g}{2}h^2\right)_\jph^-
-a_\jph^-\left(\frac{g}{2}h^2\right)_\jph^+}{a_\jph^+-a_\jph^-}\\
&+&\frac{a_\jph^+a_\jph^-}{a_\jph^+-a_\jph^-}\Big[(hu)_\jph^+-(hu)_\jph^-\Big]=\frac{g}{2}h_\jph^2,
\end{eqnarray*}
where $h_\jph:=h_\jph^+=h_\jph^-$. So, the finite volume update \eref{eq:fv-m} for the studied steady state
reads after substituting the source quadrature \eref{s1},
\begin{eqnarray*}
(\overline{hu})_j^{n+1}&=&(\overline{hu})_j^n-\frac{\dt}{\dx}\left[\frac{g}{2}(h_\jph)^2-\frac{g}{2}(h_\jmh)^2\right]+\dt\,\overline\mS_j^{(2),n} \\
&=&(\overline{hu})_j^n-\frac{\dt}{\dx}\left[\frac{g}{2}(h_\jph)^2-\frac{g}{2}(h_\jmh)^2\right]-\frac{\dt}{\dx}\, g \overline h_j({B_\jph-B_\jmh})
\\
&=&(\overline{hu})_j^n,
\end{eqnarray*}
where we have used 
\begin{equation}
\frac{(h_\jph)^2-(h_\jmh)^2}{2}=-\overline h_j^n\left(B_\jph-B_\jmh\right).
\label{eq:simpwb}
\end{equation}
It remains the verify~\eqref{eq:simpwb}. In the fully flooded cells, where $w_j>B_\jpmh$, we have
\begin{eqnarray*}
\frac{(h_\jph)^2-(h_\jmh)^2}{2}&=&\frac{h_\jph+h_\jmh}{2}\left(h_\jph-h_\jmh\right)=
\overline h_j^n\left(w_j-B_\jph-w_j+B_\jmh\right)\\
&=&-\overline h_j^n\left(B_\jph-B_\jmh\right),
\end{eqnarray*}
and thus \eref{eq:simpwb} is satisfied. In the partially flooded cells (as the one shown in Figure \ref{fig:compute_w} on
the right), $w_j<B_\jmh$, $h_\jmh=0$, and thus using \eref{eq:reccons} equation \eref{eq:simpwb} reduces to
$$
\frac{(h_\jph)^2}{2}=-\frac{\dxs_jh_\jph}{2\dx}\left(B_\jph-B_\jmh\right)=-\frac{h_\jph}{2}\dxs_j(B_x)_j,
$$
which is true since at the studied-steady situation, $\xs_j=\xs_w$, which implies that $\dxs_j=\dxs_w$, and
hence, $-\dxs_j(B_x)_j=h_\jph$.

This concludes the proof of the theorem.$\hfill\Box$

\smallskip
\begin{rem}
The draining time $\dtd$ equals the standard time step $\dt$ in all cells except at the wet/dry boundary. Therefore, the update {\rm\eref{eq:fvh}} equals the original update {\rm\eref{f1d}} almost everywhere.
\end{rem}

\smallskip
\begin{rem}
We would like to point out that the resulting scheme will clearly remain
positivity preserving if the forward Euler method in the discretization of the
ODE system {\rm\eref{f1d}} is replaced with a higher-order SSP ODE solver
(either the Runge-Kutta or the multistep one), because such solvers can be
written as a convex combination of several forward Euler steps, see \cite{GST}.
In each Runge-Kutta stage, the time step $\Delta t$ is chosen as the global time
step at the first stage. This is because the draining time $\dtd$, which is a
local cut-off to the numerical flux,  does not reduce, or even influence, the
global time step.
\end{rem}

\section{Numerical Experiments}\label{sec:exp}

Here, we set $\theta=1.3$ in the minmod function \eref{minmod}, and in
\eqref{CFL} we set $\CFL = 0.5$.

To show the effects of our new reconstruction at the boundary, we first test the numerical accuracy order using a continuous problem; then
compare our new scheme with the scheme from \cite{KP} for the oscillating lake problem and the wave run-up problem on a slopping shore. These schemes
only differ in the treatment of the dry boundary, so that the effects of the proposed modifications are highlighted. At last, we apply our
scheme to dam-break problems over a plane and a triangular hump with bottom friction. For the sake of brevity, we refer to the scheme from
\cite{KP} as KP and to our new scheme as BCKN.

Before the simulations, let us talk about the cell averages for the initial condition. Suppose that the states at cell interfaces $\mU_{j-\frac{1}{2}}$ and $\mU_{j+\frac{1}{2}}$ are given. The cell averages of momentums $(hu)_j$ are computed using the trapezoidal rule in the cells $I_j$ as
\[ (hu)_j = \frac{(hu)_{j-\frac{1}{2}} + (hu)_{j+\frac{1}{2}}}{2}.\]
As for the water height, we have to distinguish between three cases \cite{RB}. Cells $I_j$ are called wet cells if the water heights at both
cell interfaces are positive,
\[h_{j-\frac{1}{2}}>0 \quad \hbox{and} \quad  h_{j+\frac{1}{2}}>0.\]
If instead,
\[ h_{j-\frac{1}{2}}=0, \quad h_{j+\frac{1}{2}}>0 \quad \hbox{and} \quad B_{j-\frac{1}{2}}>B_{j+\frac{1}{2}},\]
we speak of cells with \textit{upward slope}.  If
\[ h_{j-\frac{1}{2}}=0, \quad h_{j+\frac{1}{2}}>0 \quad \hbox{and} \quad B_{j-\frac{1}{2}}<B_{j+\frac{1}{2}},\]
we speak of \textit{downward slope}. For the wet cells and cells with downward slope, the cell averages of water height $h_j$ are computed using the trapezoidal rule in the cells $I_j$ as
\[ h_j = \frac{h_{j-\frac{1}{2}} + h_{j+\frac{1}{2}}}{2},\]
because it is impossible to be still water states. For the adverse slope, we use the inverse function of \eref{eq:wjdry},
\[ h_j= \frac{(h_{j+\frac{1}{2}})^2}{2(B_{j+\frac{1}{2}}-B_{j-\frac{1}{2}})}, \]
to computed the cell average water height assuming the water is flat. It is easy to see from our new reconstruction \eref{eq:dxw}, \eref{eq:wjdry} and \eref{eq:dxr} can exactly reconstruct the initial still water states.

\subsection{Numerical accuracy order}\label{sec:EOC}

To compute the numerical order of accuracy of our scheme, we choose a continuous example from \cite{KP}. With computational domain $[0,1]$, the problem is subject to the gravitational constant $g=9.812$, the bottom topography
\[ B(x) = \sin ^2(\pi x), \]
the initial data
\[h(x,0)=5+e^{\cos(2\pi x)}, \quad hu(x,0)=\sin(\cos(2\pi x)),\]
and the periodic boundary conditions.
\begin{table}[htpb]
\centering
\begin{tabular}{ccccc}\hline
\quad \# points \quad & \quad $h$ error \quad & \quad EOC \quad & \quad $hu$ error \quad & \quad EOC\\\hline\hline
  25  &   5.30e-2    &         &  2.33e-1 &        \\
  50  &   1.51e-2    &   1.81  &  1.38e-1 & 0.76   \\
 100  &   4.86e-3    &   1.63  &  4.43e-2 & 1.64   \\
 200  &   1.40e-3    &   1.80  &  1.14e-2 & 1.95   \\
 400  &   3.59e-4    &   1.96  &  2.84e-3 & 2.01   \\
 800  &   8.93e-5    &   2.01  &  7.05e-4 & 2.01  \\\hline
\end{tabular}
\caption{Accuracy checking: Experimental order of convergence(EOC) measured in
the $L^1$-norm.}
\label{tab:accuracy}
\end{table}

 The reference solution is computed on a grid with $12800$ cells. The numerical result is shown in the Table \ref{tab:accuracy} at time $t=0.1$. The result  confirm that our scheme is second-order accurate.

\subsection{Still and oscillating lakes }\label{sec:ol}

In this section, we consider present a test case proposed in \cite{ABBKP}. It describes the situation where the
``lake at rest'' \eref{sss} and ``dry lake'' \eref{dss} are combined in the domain $[0,1]$ with the bottom
topography given by
\begin{equation}
B(x)=\frac{1}{4}-\frac{1}{4}\cos((2x-1)\pi),
\label{5.1}
\end{equation}
and the following initial data:
\begin{equation}
h(x,0)=\max\left(0,0.4-B(x)\right),\qquad u(x,0)\equiv 0.
\label{5.2}
\end{equation}

\begin{figure}[htpb]
\centering
\includegraphics[width=0.8\linewidth]{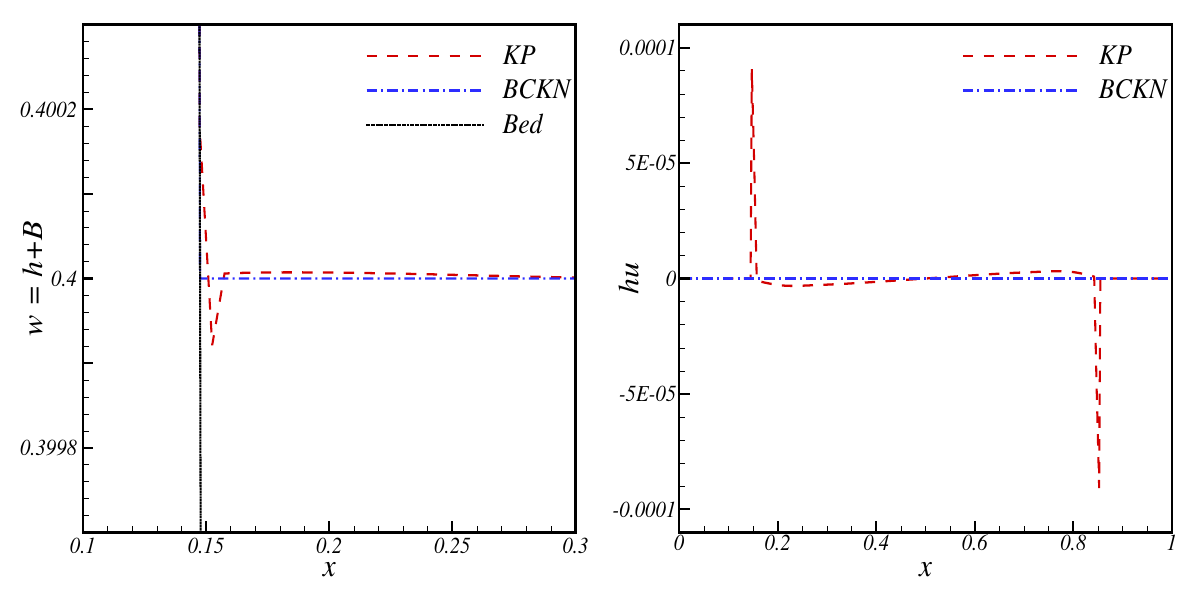}
\caption{Lake at rest. Left: free surface $h+B$; Right: Discharge $hu$ (cf. Table \ref{tab:err_lar}).}
\label{fig:lar_cr}
\end{figure}

We compute the numerical solution by the KP and BCKN schemes with $200$ points at the final time $T=19.87$. The
results are shown in Figure \ref{fig:lar_cr} and Table \ref{tab:err_lar}.
As one can clearly see there, the KP scheme introduces some oscillations at the boundary, whereas the BCKN scheme is perfectly well-balanced which means that our new initial
data reconstruction method can exactly preserve the well-balanceed property not only in the wet region but also in the dry region. And the influence on the solutions away from the wet/dry front is also visible  because of oscillations at the boundary produced by the KP schme.
\begin{table}[htpb]
\centering
\begin{tabular}{ccc}\hline
\quad scheme \quad & \quad $L^\infty$ error of $h$ \quad & \quad $L^\infty$ error of $hu$\\\hline\hline
KP  & 7.88e-5   & 9.08e-5 \\
BCKN & 3.33e-16   & 5.43e-16\\
\hline
\end{tabular}
\caption{Errors in the computation of the steady state (cf. Figure \ref{fig:lar_cr})}
\label{tab:err_lar}
\end{table}

We now consider a sinusoidal perturbation of the steady state \eref{5.1}, \eref{5.2} by taking
\begin{equation*}
h(x,0)=\max\left(0,0.4+\frac{\sin\left(4x-2-\max(0,-0.4+B(x))\right)}{25}-B(x)\right).
\end{equation*}
As in \cite{ABBKP}, we set the final time to be $T=19.87$. At this time, the wave has its maximal height at the
left shore after some oscillations.

\begin{figure}[htpb]
\centering
\includegraphics[width=0.8\linewidth]{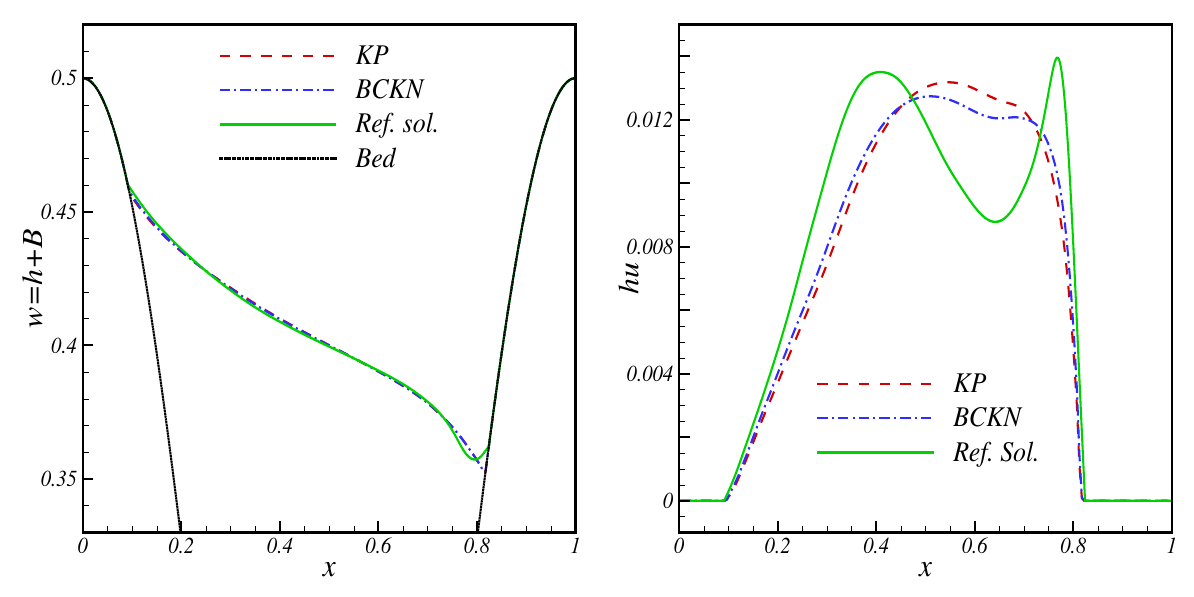}
\caption{Oscillating lake. Left: Free surface $h+B$; Right: Discharge $hu$. Comparison of KP and BCKN schemes with the reference solution.}
\label{fig:ol}
\end{figure}
\begin{figure}[htpb]
\centering
\includegraphics[width=0.8\linewidth]{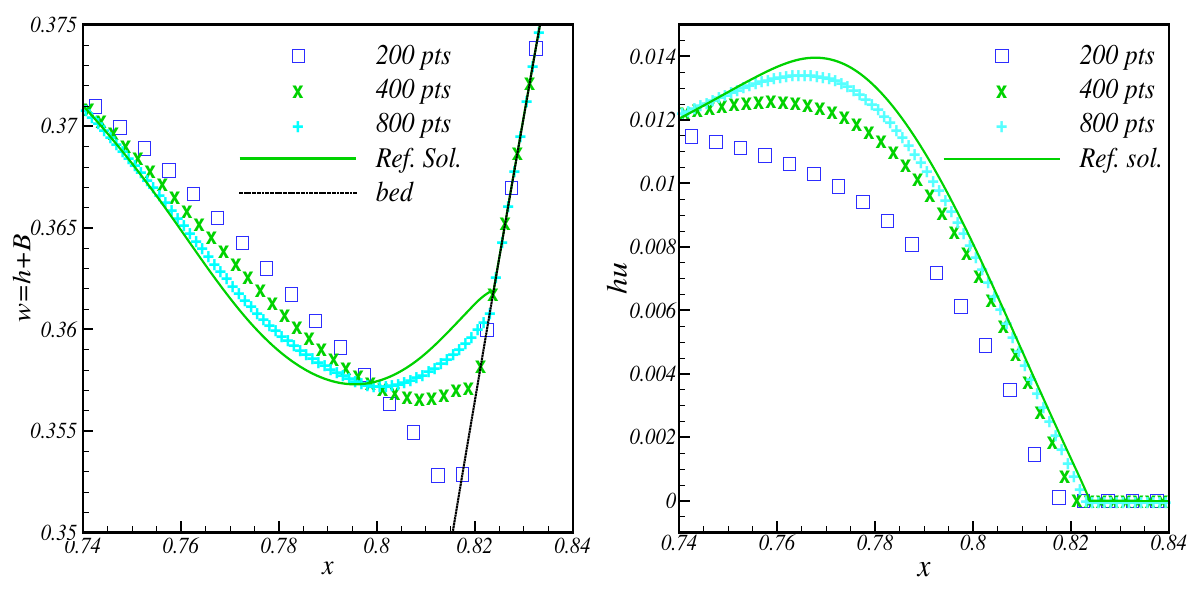}
\caption{Oscillating lake, zoom at the right wet/dry front. BCKN solutions with 200, 400, 800 points and reference solution ($12800$ points). Left: Free surface $h+B$; Right: Discharge $hu$.}
\label{fig:o2}
\end{figure}

In Figure \ref{fig:ol} we compare the results obtained by the BCKN and KP schemes with $200$ points with a reference solution (computed using $12800$ points). Table \ref{tab:conv_new} shows the experimental accuracy order for the two different schemes. One can clearly see that both KP and BCKN scheme can produce good results and acceptable numerical order. In Figure \ref{fig:o2} we show a zoom of BCKN solutions for $x\in[0.74,0.84]$ with 200, 400 and 800 points, which converge nicely to the reference solution. In particular, the discharge converges without any oscillations.
\begin{table}[htpb]
\centering
\begin{tabular}{ccccc}\hline
\# points&$h$ error&EOC&$hu$ error&EOC\\\hline\hline
        25 &  9.48e-3 &            &  1.47e-2 &           \\
        50 &  2.81e-3 &       1.75 &  7.26e-3 &       1.02 \\
       100 &  1.65e-3 &       0.77 &  2.46e-3 &       1.56 \\
       200 &  7.88e-4 &       1.06 &  1.59e-3 &       0.63 \\
       400 &  3.33e-4 &       1.24 &  6.19e-4 &       1.36 \\
       800 &  1.26e-4 &       1.40 &  2.27e-4 &       1.45 \\
 \hline
 KP scheme\\
 \hline
        25 &  7.55e-3 &          &  1.31e-2 &           \\
        50 &  2.27e-3 &     1.74 &  6.04e-3 &       1.11 \\
       100 &  1.45e-3 &     0.65 &  2.35e-3 &       1.36 \\
       200 &  6.77e-4 &     1.09 &  1.31e-3 &       0.84 \\
       400 &  2.71e-4 &     1.32 &  5.04e-4 &       1.38 \\
       800 &  1.04e-4 &     1.38 &  1.87e-4 &       1.43 \\
\hline
BCKN scheme\\
 \hline
 \end{tabular}
\caption{Oscillating lake: Experimental order of convergence  measured in
the $L^1$-norm.}
\label{tab:conv_new}
\end{table}

\subsection{Wave run-up on a sloping shore}\label{sec:shore}

This test describes the run-up and reflection of a wave on a mounting slope. It was proposed in \cite{Sy} and reference solutions can be found, for example, in \cite{BNL,RB,Synolakis1986}.

The initial data are
$$
H_0(x)=\max\left\{D+\delta\,{\rm sech}^2(\gamma(x-x_a)),B(x)\right\},\qquad u_0(x)=\sqrt{\frac{g}{D}}\,H_0(x),
$$
and the bottom topography is
$$
B(x)=\begin{cases}0,& \text{if~$x<2x_a$},\\\dfrac{x-2x_a}{19.85},&\text{otherwise}.\end{cases}
$$
As in \cite{BNL,RB}, we set
$$
D=1,\quad\delta=0.019,\quad\gamma=\sqrt{\frac{3\delta}{4D}},\quad
x_a=\sqrt{\frac{4D}{3\delta}}\,{\rm arccosh}\big(\sqrt{20}\big).
$$
The computational domain is $[0,80]$ and the number of grid cells is $200$.

\begin{figure}[htpb]
\centering
\includegraphics[width=0.8\linewidth]{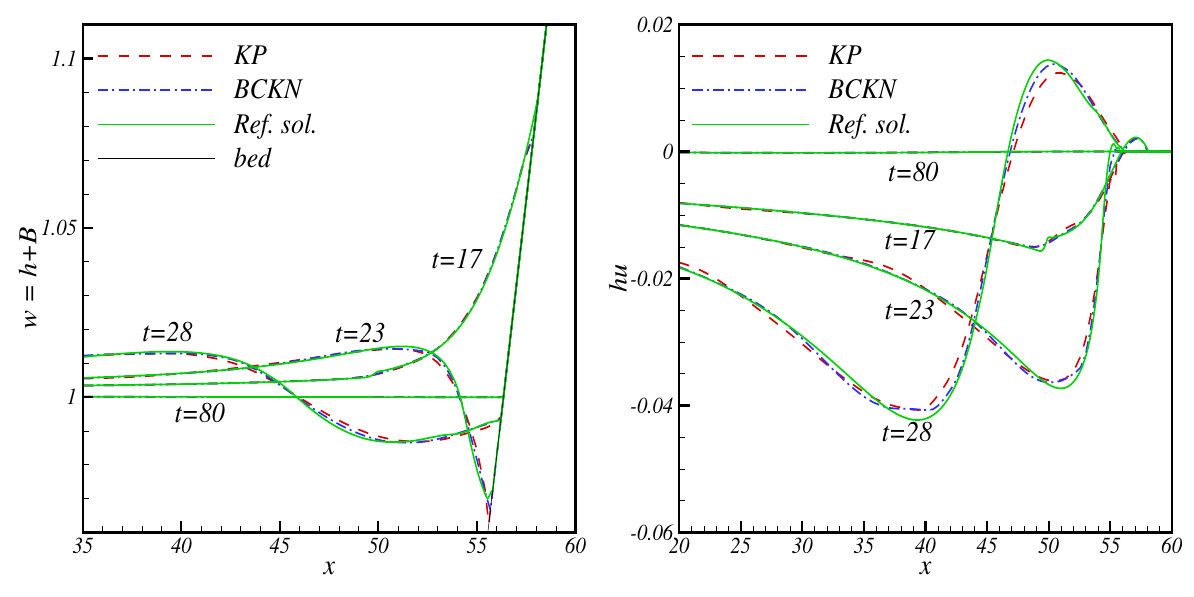}
\caption{Wave run-up on a sloping shore. KP, BCKN and reference solutions at times 17, 23, 28 and 80. Left: free surface $w=h+B$; Right: discharge $hu$.}
\label{fig:shore-dynamic}
\end{figure}

\begin{figure}[htpb]
\centering
\includegraphics[width=0.8\linewidth]{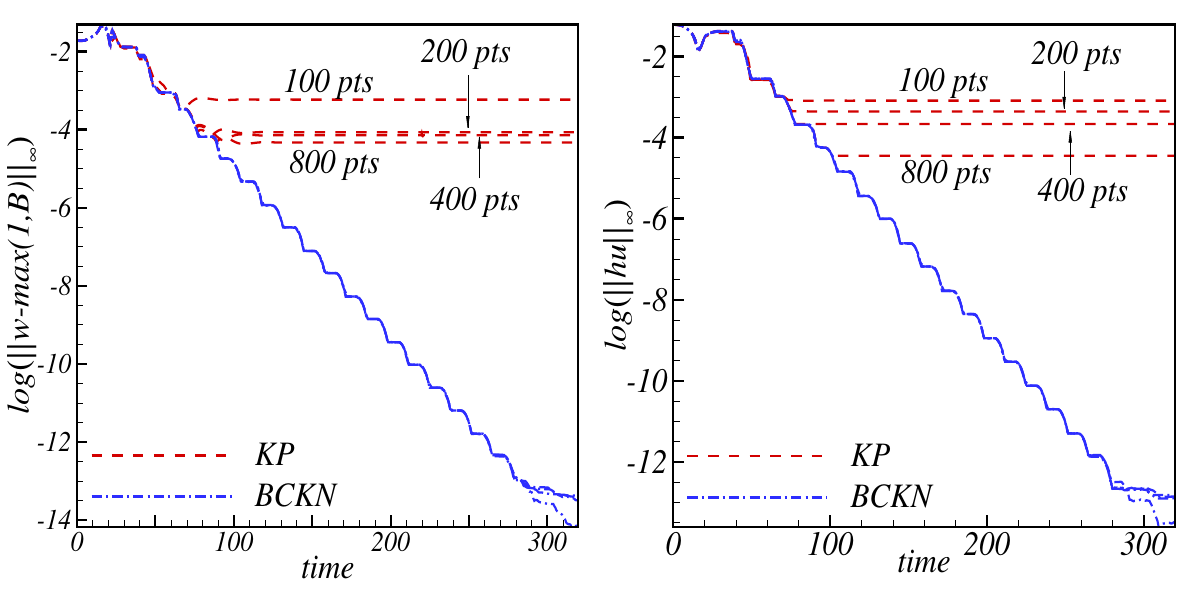}
\caption{Wave run-up on a sloping shore: deviation from stationary state. Left: free surface $\log(||w-\max(1,B)||_{\infty})$; Right: discharge $\ln(||hu||_{\infty})$. KP scheme (dashed) and BCKN scheme (dash-dot). Long time convergence of KP scheme stalls.}
\label{fig:shore-long-time}
\end{figure}

Figure \ref{fig:shore-dynamic} shows the free surface and discharge
computed by both BCNK and KP schemes for different times. The reference solution
is computed using 2000 points. A wave is running up the shore at time $t=17$, and running down at $t=23$. At time $t=80$ a steady state is reached. In the dynamic phase (up to time $t=28$), both schemes provide satisfactory solutions.
In Figure \ref{fig:shore-long-time} we study the long time decay towards equilibrium for different grid size resolutions. While the BCKN solutions decay up to machine accuracy, the long time convergence of the KP scheme comes to a halt. A brief check reveals that the deviation from equilibrium is roughly of the size of the truncation error of the KP scheme.

\subsection{Dam-break over a plane}\label{sec:DB_plane}

Here we study three dam breaks over inclined planes with various inclination angles. These test cases have been
previously considered in \cite{GPC2007, xing2010}.

The domain is $[-15,15]$, the bottom topography is given by
	 	\[B(x)=-x\tan \alpha \]
where $\alpha$ is the inclination angle. The initial data are
	 	\[u(x,0)=0,\ \ \ \ \ \ \ \ h(x,0)=\left\{ \begin{matrix}
   1-B(x), & x<0,  \\
   0, & \text{otherwise}\text{.}  \\
\end{matrix} \right.\]

At $x=15$ we impose a free flow boundary condition, and at $x=-15$ we set the discharge to zero. The plane is either flat ($\alpha =0$), inclined uphill ($\alpha =\pi/60$), or downhill ($\alpha =-\pi/60$).

\begin{figure}[htpb]
\centering
\includegraphics[width=0.8\linewidth]{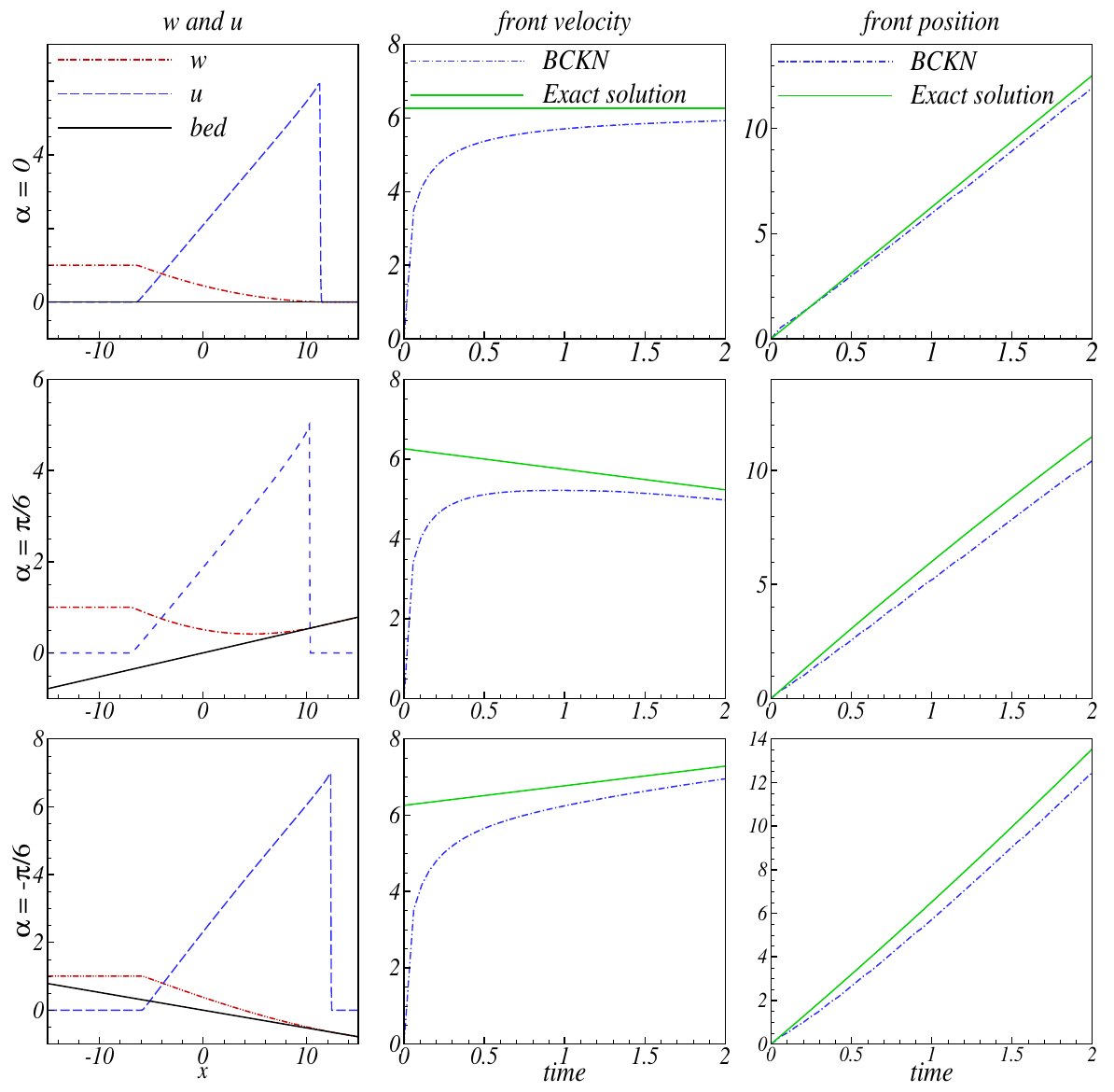}
\caption{Dam-break over a plane. Left : the numerical solution of $w=h+B$ and $u$; Middle: the front position; Right: the front velocity.}
\label{fig:plane}
\end{figure}

We run the simulation until time $t =2$, with 200 uniform cells. The numerical results are displayed in Figure \ref{fig:plane}, for inclination angles  $\alpha =0$, $\pi/60$ and $-\pi/60$, from top to bottom. The left column shows $h$ and $u$, the central column the front position and the right column the front velocity. We also display the exact front positions and velocities (see \cite{GPC2007}) given by
\[{{x}_{f}}(t)=2t\sqrt{g\cos (\alpha )}-\frac{1}{2}g{{t}^{2}}\tan (\alpha ), \quad {{u}_{f}}(t)=2\sqrt{g\cos (\alpha )}-gt\tan (\alpha ).\]
As suggested in \cite{xing2010}, we define the numerical front position to be the first cell (counted from right to left) where the water height exceeds $\epsilon=10^{-9}$.
While the BCKN scheme, which is only second order accurate, cannot fully match the resolution of the third and fifth order schemes in \cite{GPC2007, xing2010}, it still performs reasonably well. What we would like to stress here is that the new scheme, which was designed to be well balanced near wet/dry equilibrium states, is also robust for shocks running into dry areas.

\subsection{Laboratory dam-break over a triangular hump}\label{sec:hump}

\begin{figure}[htpb]
\centering
\includegraphics[width=0.8\linewidth]{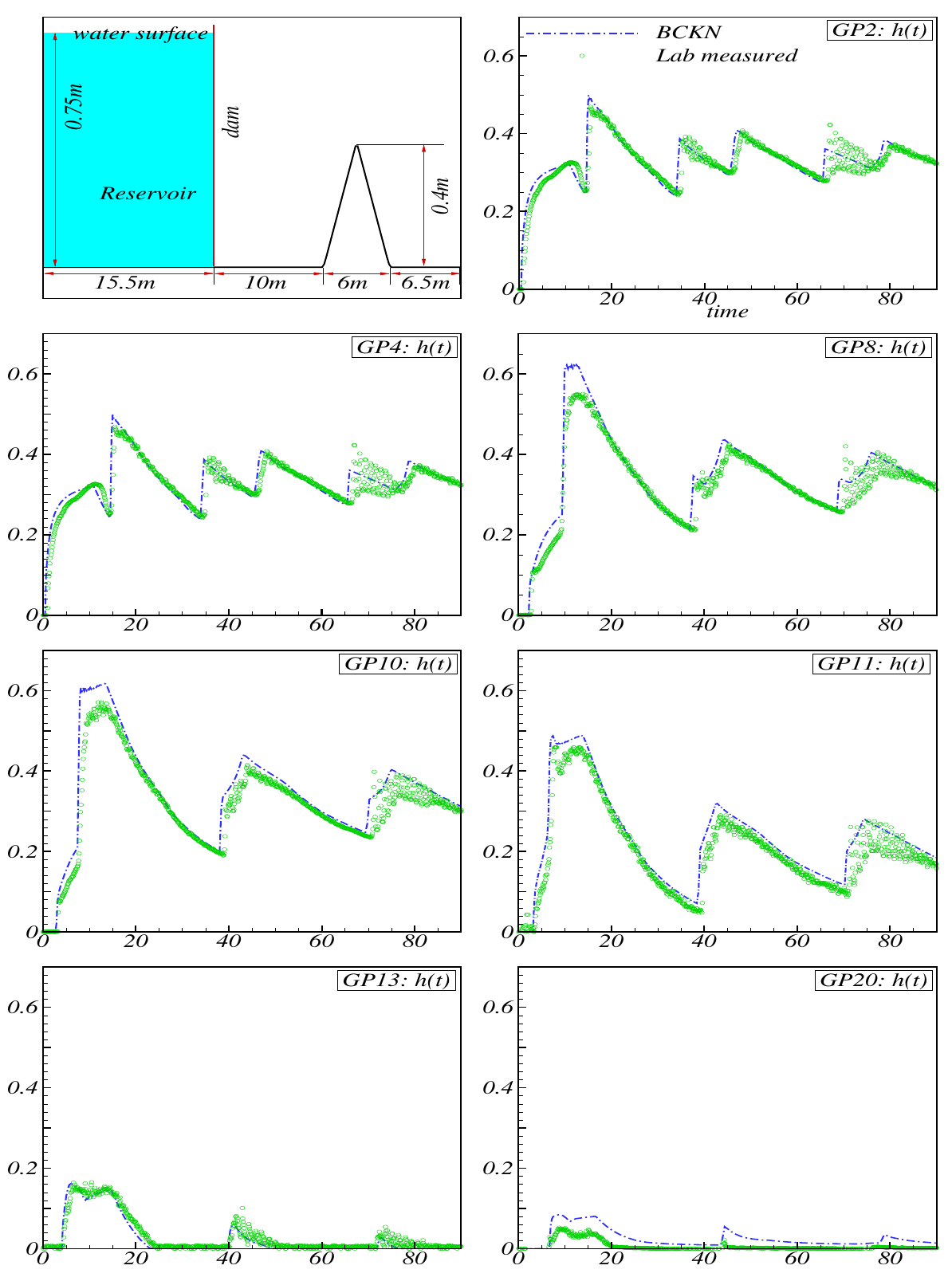}
\caption{Laboratory dam-break inundation over dry bed: experimental setup and the comparison of simulated and observed water depth versus time at 7 gauge points.}
\label{fig:hump}
\end{figure}

We apply our scheme to a laboratory test of a dam-break inundation over a
triangular hump which is recommended by the Europe, the Concerted Action on
Dam-Break Modeling (CADAM) project \cite{Morris2000}. The problem consider the
friction effect and then the corresponding governing equation (\ref{1d}) is
changed to be
\begin{eqnarray}
\left\{\begin{array}{l}\displaystyle{h_t+(hu)_x=0,}\\[0.3ex]
\displaystyle{(hu)_t +\Big(hu^2+\frac{1}{2}gh^2\Big)_x=-ghB_x - \tau_b/\rho },\end{array}\right.
\label{eqn:friction}
\end{eqnarray}
where $\tau_b =  \rho c_f u|u|$ represents the energy dissipation effect and are estimated from bed roughness on the flow, $\rho$ is the density of water and $c_f=gn^2/h^{1/3}$ represents the bed roughness coefficient with $n$ being the Manning coefficient. For small water depths, the bed friction term dominates the other terms in the momentum equation, due to the presence of $h^{1/3}$ in the denominator. To simplify the update of the momentum, we first update the solution using our new positivity preserving and well-balanced scheme stated in section 2, 3 and 4 without the bed friction effect, and then retain the local acceleration from the only bed friction terms.
\begin{equation}
(hu)_t = - \tau_b/\rho = -c_f u|u|=- \frac{gn^2u|u|}{h^{1/3}}.
\end{equation}
A partially implicit approach \cite{Liang2009,Singh2011} is used for the discretization of the above equation as
\begin{equation}
\frac{(hu)^{n + 1} -(\tilde{h}\tilde{u})^{n + 1} }{\Delta t}  = - \frac{gn^2(hu)^{n+1}|\tilde{u}^{n+1}|}{(\tilde{h}^{n+1})^{4/3}}.
\end{equation}
Resolving this for $(hu)^{n + 1}$, we obtain 
\begin{equation}
(hu)^{n + 1} =\frac{(\tilde{h}\tilde{u})^{n + 1}}{1+\Delta t gn^2|\tilde{u}^{n+1}|/(\tilde{h}^{n+1})^{4/3}}= \frac{(\tilde{h}\tilde{u})^{n + 1} (\tilde{h}^{n+1})^{4/3}}{(\tilde{h}^{n+1})^{4/3}+\Delta t gn^2|\tilde{u}^{n+1}|},
\end{equation}
where $\tilde{h}$ and $\tilde{u}$ are given using our above stated scheme without friction term.  The initial conditions and geometry (Figure \ref{fig:hump}) were identical to those used by \cite{Liang2009,Singh2011}. The experiment was conducted in a 38-m-long channel. The dam was located at 15.5 m, with a still water surface of 0.75 m in the reservoir. A symmetric triangular obstacle 6.0 m long and 0.4 m high was installed 13.0 m downstream of the dam. The floodplain was fixed and initially dry, with reflecting boundaries and a free outlet. The Manning coefficient $n$ was 0.0125, adopted from \cite{Liang2009}. The flow depth was measured at seven stations, GP2, GP4, GP8, GP10, GP11, GP13, and GP20, respectively, located at 2, 4, 8, 10, 11, 13, and 20 m downstream of the dam, as shown in Figure  \ref{fig:hump}. The simulation was conducted for 90 seconds.

The numerical predictions using 200 points are shown in Figure \ref{fig:hump}. The comparison between the numerical results and measurements
is satisfactory at all gauge points and the wet/dry transitions are resolved sharply (compare with \cite{Liang2009,Singh2011} and
the references therein). This confirms the effectiveness of the current scheme together with the implicit method for discretization of the
friction term, even near wet/dry fronts.

\section{Conclusion}

In this paper, we designed a special reconstruction of the water level at
wet/dry fronts, in the framework of the second-order semi-discrete
central-upwind scheme and a continuous, piecewise linear discretisation of the
bottom topography. The proposed reconstruction is conservative, well-balanced
and positivity preserving for both wet and dry cells. The positivity of the
computed water height is ensured by cutting the outflux across partially flooded edges
at the draining time, when the cell has run empty. Several numerical examples
demonstrate the experimental order of convergence and the well-balancing
property of the new scheme, and we also show a case where the prerunner of the
scheme fails to converge to equilibrium. The new scheme is robust for shocks
running into dry areas and for simulations including Manning's bottom friction
term, which is singular at the wet/dry front.

\medskip
{\bf Acknowledgment.} The first ideas for this work were discussed by the authors at a meeting at the ``Mathematisches Forschungsinstitut
Oberwolfach''. The authors are grateful for the support and inspiring atmosphere there. The research of A. Kurganov was supported in part by
the NSF Grant DMS-1115718 and the ONR Grant N000141210833. The research of A. Bollermann, G. Chen and S. Noelle was supported by DFG Grant
NO361/3-1 and No361/3-2. G. Chen is partially supported by the National Natural Science Foundation of China (No. 11001211, 51178359).

\bibliographystyle{amsplain}
\bibliography{Xbib}

\providecommand{\bysame}{\leavevmode\hbox to3em{\hrulefill}\thinspace}
\providecommand{\MR}{\relax\ifhmode\unskip\space\fi MR }
\providecommand{\MRhref}[2]{%
  \href{http://www.ams.org/mathscinet-getitem?mr=#1}{#2}
}
\providecommand{\href}[2]{#2}
\begin{thebibliography}{10}

\bibitem{ABBKP}
E.~Audusse, F.~Bouchut, M.-O. Bristeau, R.~Klein, and B.~Perthame, \emph{A fast
  and stable well-balanced scheme with hydrostatic reconstruction for shallow
  water flows}, SIAM J. Sci. Comput. \textbf{25} (2004), 2050--2065.

\bibitem{BNL}
A.~Bollermann, S.~Noelle, and M.~Luk\'{a}\v{c}ov\'{a}-Medvid'ov\'{a},
  \emph{Finite volume evolution {G}alerkin methods for the shallow water
  equations with dry beds}, Commun. Comput. Phys. \textbf{10} (2011), 371--404.

\bibitem{SV}
A.J.C. de~Saint-Venant, \emph{\em {T}h\'{e}orie du mouvement non-permanent des
  eaux, avec application aux crues des rivi\`{e}re at \`{a} l'introduction des
  war\'{e}es dans leur lit}, C.R. Acad. Sci. Paris \textbf{73} (1871),
  147--154.

\bibitem{GPC2007}
J.M. Gallardo, C.~Par\'{e}s, and M.~Castro, \emph{On a well-balanced high-order
  finite volume scheme for shallow water equations with topography and dry
  areas}, J. Comput. Phys. \textbf{2227} (2007), 574--601.

\bibitem{GHS}
T.~Gallou\"et, J.-M. H\'erard, and N.~Seguin, \emph{Some approximate {G}odunov
  schemes to compute shallow-water equations with topography}, Comput. \&
  Fluids \textbf{32} (2003), 479--513.

\bibitem{GR2}
E.~Godlewski and P.-A. Raviart, \emph{Numerical approximation of hyperbolic
  systems of conservation laws}, Springer-Verlag, New York, 1996.

\bibitem{GST}
S.~Gottlieb, C.-W. Shu, and E.~Tadmor, \emph{High order time discretization
  methods with the strong stability property}, SIAM Review \textbf{43} (2001),
  89--112.

\bibitem{Jin}
S.~Jin, \emph{A steady-state capturing method for hyperbolic system with
  geometrical source terms}, M2AN Math. Model. Numer. Anal. \textbf{35} (2001),
  631--645.

\bibitem{JW}
S.~Jin and X.~Wen, \emph{Two interface-type numerical methods for computing
  hyperbolic systems with geometrical source terms having concentrations}, SIAM
  J. Sci. Comput. \textbf{26} (2005), 2079--2101.

\bibitem{Kr}
D.~Kr\"oner, \emph{Numerical schemes for conservation laws}, Wiley, Chichester,
  1997.

\bibitem{KL}
A.~Kurganov and D.~Levy, \emph{Central-upwind schemes for the {S}aint-{V}enant
  system}, M2AN Math. Model. Numer. Anal. \textbf{36} (2002), 397--425.

\bibitem{KP}
A.~Kurganov and G.~Petrova, \emph{A second-order well-balanced positivity
  preserving central-upwind scheme for the {S}aint-{V}enant system}, Commun.
  Math. Sci. \textbf{5} (2007), no.~1, 133--160.

\bibitem{LeVbook}
R.~LeVeque, \emph{Finite volume methods for hyperbolic problems}, Cambridge
  Texts in Applied Mathematics,Cambridge University Press, 2002.

\bibitem{LeV}
R.J. LeVeque, \emph{Balancing source terms and flux gradients in
  high-resolution {G}odunov methods: the quasi-steady wave-propagation
  algorithm}, J. Comput. Phys. \textbf{146} (1998), 346--365.

\bibitem{Liang2009}
Q.~Liang and F.~Marche, \emph{Numerical resolution of well-balanced shallow
  water equations with complex source terms}, Adv. Water Resour. \textbf{32}
  (2009), no.~6, 873--884.

\bibitem{LieNoe2}
K.-A. Lie and S.Noelle, \emph{On the artificial compression method for
  second-order nonoscillatory central difference schemes for systems of
  conservation laws}, SIAM J. Sci. Comput. \textbf{24} (2003), 1157--1174.

\bibitem{Morris2000}
M.~Morris, \emph{\emph{{C}{A}{D}{A}{M}:} concerted action on dambreak modeling
  - final report}, HR Wallingford, 2000.

\bibitem{NT}
H.~Nessyahu and E.~Tadmor, \emph{Non-oscillatory central differencing for
  hyperbolic conservation laws}, J. Comput. Phys. \textbf{87} (1990), 408--463.

\bibitem{NPPN}
S.~Noelle, N.~Pankratz, G.~Puppo, and J.~Natvig, \emph{Well-balanced finite
  volume schemes of arbitrary order of accuracy for shallow water flows}, J.
  Comput. Phys. \textbf{213} (2006), 474--499.

\bibitem{PS}
B.~Perthame and C.~Simeoni, \emph{A kinetic scheme for the {S}aint-{V}enant
  system with a source term}, Calcolo \textbf{38} (2001), 201--231.

\bibitem{RB}
M.~Ricchiuto and A.~Bollermann, \emph{Stabilized residual distribution for
  shallow water simulations}, J. Comput. Phys. \textbf{228} (2009), 1071--1115.

\bibitem{Rus1}
G.~Russo, \emph{Central schemes for balance laws}, Hyperbolic problems: theory,
  numerics, applications, Vol. I, II (Magdeburg, 2000), 821-829, Internat. Ser.
  Numer. Math., 140, 141, Birkh\"auser, Basel, 2001.

\bibitem{Rus2}
\bysame, \emph{Central schemes for conservation laws with application to
  shallow water equations}, in Trends and applications of mathematics to
  mechanics: STAMM 2002, S. Rionero and G. Romano (eds.), 225--246,
  Springer-Verlag Italia SRL, 2005.

\bibitem{Singh2011}
J.~Singh, M.S. Altinakar, and Y.~Ding, \emph{Two-dimensional numerical modeling
  of dam-break flows over natural terrain using a central explicit scheme},
  Adv. Water Resour. \textbf{34} (2011), 1366--1375.

\bibitem{Swe}
P.~K. Sweby, \emph{High resolution schemes using flux limiters for hyperbolic
  conservation laws}, SIAM J. Numer. Anal. \textbf{21} (1984), 995--1011.

\bibitem{Synolakis1986}
C.~E. Synolakis., \emph{The runup of long waves}, Ph.D. thesis, California
  Institute of Technology, 1986.

\bibitem{Sy}
C.~E. Synolakis, \emph{The runup of solitary waves}, J. Fluid Mech.
  \textbf{185} (1987), 523--545.

\bibitem{TNGH}
Y.C. Tai, S.~Noelle, J.M.N.T. Gray, and K.~Hutter, \emph{Shock-capturing and
  front-tracking methods for granular avalanches}, J. Comput. Phys.
  \textbf{175} (2002), 269--301.

\bibitem{vLeV}
B.~van Leer, \emph{Towards the ultimate conservative difference scheme, {V}.
  {A} second order sequel to {G}odunov's method}, J. Comput. Phys. \textbf{32}
  (1979), 101--136.

\bibitem{XS1}
Y.~Xing and C.-W. Shu, \emph{High order finite difference {WENO} schemes with
  the exact conservation property for the shallow water equations}, J. Comput.
  Phys. \textbf{208} (2005), 206--227.

\bibitem{XS2}
\bysame, \emph{A new approach of high order well-balanced finite volume {WENO}
  schemes and discontinuous {G}alerkin methods for a class of hyperbolic
  systems with source terms}, Commun. Comput. Phys. \textbf{1} (2006),
  100--134.

\bibitem{xing2010}
Y.~Xing, X.~Zhang, and C.W. Shu, \emph{Positivity-preserving high order
  well-balanced discontinuous {G}alerkin methods for the shallow water
  equations}, Adv. Water Resour. \textbf{33} (2010), no.~12, 1476--1493.

\end{thebibliography}

\end{spacing}

\end{document}